\documentclass[10.5pt]{article}
\usepackage{amsmath}
\usepackage{amssymb}
\usepackage{amscd}
\usepackage{xspace}
\usepackage{verbatim}
\newcommand{\mysection}[1]{
\section{#1}\setcounter{equation}{0}}
\title{\bf Separable solutions of some quasilinear equations with source reaction}
\author{{\bf Marie-Fran\c{c}oise Bidaut-V\'eron}\\
{\small Department of Mathematics, Universit\'e Fran\c{c}ois Rabelais,  Tours,  FRANCE}\\[2mm]
{\bf Mustapha Jazar\footnote{Supported by a grant from the Lebanese University}}\\
{\small Department of Mathematics, Universit\'e Libanaise, Beyrouth, LIBAN}\\[2mm]
 {\bf Laurent V\'eron}\\
{\small Department of Mathematics, Universit\'e Fran\c{c}ois Rabelais,  Tours,  FRANCE}
}

\date{}
\begin{document}
\maketitle
\noindent {\small {\bf Abstract} We study the existence of singular  solutions to the equation
$-div (|Du|^{p-2}D u)=| u|^{q-1}u$ under
the form $u(r,\theta)=r^{-\beta}\omega(\theta)$, $r>0,\theta\in S^{N-1}$. We prove the existence of an exponent $q$ below which no positive solutions can exist. If the dimension is $2$ we use a dynamical system approach to construct solutions.} \smallskip

\vspace{1mm}
\hspace{.05in}

\newcommand{\txt}[1]{\;\text{ #1 }\;}
\newcommand{\tbf}{\textbf}
\newcommand{\tit}{\textit}
\newcommand{\tsc}{\textsc}
\newcommand{\trm}{\textrm}
\newcommand{\mbf}{\mathbf}
\newcommand{\mrm}{\mathrm}
\newcommand{\bsym}{\boldsymbol}
\newcommand{\scs}{\scriptstyle}
\newcommand{\sss}{\scriptscriptstyle}
\newcommand{\txts}{\textstyle}
\newcommand{\dsps}{\displaystyle}
\newcommand{\fnz}{\footnotesize}
\newcommand{\scz}{\scriptsize}
\newcommand{\be}{
\begin{equation}
}
\newcommand{\bel}[1]{
\begin{equation}
\label{#1}}
\newcommand{\ee}{
\end{equation}
}
\newcommand{\eqnl}[2]{
\begin{equation}
\label{#1}{#2}
\end{equation}
}
\newtheorem{subn}{\name}
\renewcommand{\thesubn}{}
\newcommand{\bsn}[1]{\def\name{#1}
\begin{subn}}
\newcommand{\esn}{
\end{subn}}
\newtheorem{sub}{\name}[section]
\newcommand{\dn}[1]{\def\name{#1}}   
\newcommand{\bs}{
\begin{sub}}
\newcommand{\es}{
\end{sub}}
\newcommand{\bsl}[1]{
\begin{sub}\label{#1}}
\newcommand{\bth}[1]{\def\name{Theorem}
\begin{sub}\label{t:#1}}
\newcommand{\blemma}[1]{\def\name{Lemma}
\begin{sub}\label{l:#1}}
\newcommand{\bcor}[1]{\def\name{Corollary}
\begin{sub}\label{c:#1}}
\newcommand{\bdef}[1]{\def\name{Definition}
\begin{sub}\label{d:#1}}
\newcommand{\bprop}[1]{\def\name{Proposition}
\begin{sub}\label{p:#1}}
\newcommand{\R}{\eqref}
\newcommand{\rth}[1]{Theorem~\ref{t:#1}}
\newcommand{\rlemma}[1]{Lemma~\ref{l:#1}}
\newcommand{\rcor}[1]{Corollary~\ref{c:#1}}
\newcommand{\rdef}[1]{Definition~\ref{d:#1}}
\newcommand{\rprop}[1]{Proposition~\ref{p:#1}}
\newcommand{\BA}{
\begin{array}}
\newcommand{\EA}{
\end{array}}
\newcommand{\BAN}{\renewcommand{\arraystretch}{1.2}
\setlength{\arraycolsep}{2pt}
\begin{array}}
\newcommand{\BAV}[2]{\renewcommand{\arraystretch}{#1}
\setlength{\arraycolsep}{#2}
\begin{array}}
\newcommand{\BSA}{
\begin{subarray}}
\newcommand{\ESA}{
\end{subarray}}
\newcommand{\BAL}{
\begin{aligned}}
\newcommand{\EAL}{
\end{aligned}}
\newcommand{\BALG}{
\begin{alignat}}
\newcommand{\EALG}{
\end{alignat}}
\newcommand{\BALGN}{
\begin{alignat*}}
\newcommand{\EALGN}{
\end{alignat*}}
\newcommand{\note}[1]{\textit{#1.}\hspace{2mm}}
\newcommand{\Proof}{\note{Proof}}
\newcommand{\qeda}{\hspace{10mm}\hfill $\square$}
\newcommand{\qed}{\\
${}$ \hfill $\square$}
\newcommand{\Remark}{\note{Remark}}
\newcommand{\modin}{$\,$\\
[-4mm] \indent}
\newcommand{\forevery}{\quad \forall}
\newcommand{\set}[1]{\{#1\}}
\newcommand{\setdef}[2]{\{\,#1:\,#2\,\}}
\newcommand{\setm}[2]{\{\,#1\mid #2\,\}}
\newcommand{\lra}{\longrightarrow}
\newcommand{\lla}{\longleftarrow}
\newcommand{\llra}{\longleftrightarrow}
\newcommand{\Lra}{\Longrightarrow}
\newcommand{\Lla}{\Longleftarrow}
\newcommand{\Llra}{\Longleftrightarrow}
\newcommand{\warrow}{\rightharpoonup}
\newcommand{
\paran}[1]{\left (#1 \right )}
\newcommand{\sqbr}[1]{\left [#1 \right ]}
\newcommand{\curlybr}[1]{\left \{#1 \right \}}
\newcommand{\abs}[1]{\left |#1\right |}
\newcommand{\norm}[1]{\left \|#1\right \|}
\newcommand{
\paranb}[1]{\big (#1 \big )}
\newcommand{\lsqbrb}[1]{\big [#1 \big ]}
\newcommand{\lcurlybrb}[1]{\big \{#1 \big \}}
\newcommand{\absb}[1]{\big |#1\big |}
\newcommand{\normb}[1]{\big \|#1\big \|}
\newcommand{
\paranB}[1]{\Big (#1 \Big )}
\newcommand{\absB}[1]{\Big |#1\Big |}
\newcommand{\normB}[1]{\Big \|#1\Big \|}

\newcommand{\thkl}{\rule[-.5mm]{.3mm}{3mm}}
\newcommand{\thknorm}[1]{\thkl #1 \thkl\,}
\newcommand{\trinorm}[1]{|\!|\!| #1 |\!|\!|\,}
\newcommand{\bang}[1]{\langle #1 \rangle}
\def\angb<#1>{\langle #1 \rangle}
\newcommand{\vstrut}[1]{\rule{0mm}{#1}}
\newcommand{\rec}[1]{\frac{1}{#1}}
\newcommand{\opname}[1]{\mbox{\rm #1}\,}
\newcommand{\supp}{\opname{supp}}
\newcommand{\dist}{\opname{dist}}
\newcommand{\myfrac}[2]{{\displaystyle \frac{#1}{#2} }}
\newcommand{\myint}[2]{{\displaystyle \int_{#1}^{#2}}}
\newcommand{\mysum}[2]{{\displaystyle \sum_{#1}^{#2}}}
\newcommand {\dint}{{\displaystyle \int\!\!\int}}
\newcommand{\q}{\quad}
\newcommand{\qq}{\qquad}
\newcommand{\hsp}[1]{\hspace{#1mm}}
\newcommand{\vsp}[1]{\vspace{#1mm}}
\newcommand{\ity}{\infty}
\newcommand{\prt}{
\partial}
\newcommand{\sms}{\setminus}
\newcommand{\ems}{\emptyset}
\newcommand{\ti}{\times}
\newcommand{\pr}{^\prime}
\newcommand{\ppr}{^{\prime\prime}}
\newcommand{\tl}{\tilde}
\newcommand{\sbs}{\subset}
\newcommand{\sbeq}{\subseteq}
\newcommand{\nind}{\noindent}
\newcommand{\ind}{\indent}
\newcommand{\ovl}{\overline}
\newcommand{\unl}{\underline}
\newcommand{\nin}{\not\in}
\newcommand{\pfrac}[2]{\genfrac{(}{)}{}{}{#1}{#2}}

\def\ga{\alpha}     \def\gb{\beta}       \def\gg{\gamma}
\def\gc{\chi}       \def\gd{\delta}      \def\ge{\epsilon}
\def\gth{\theta}                         \def\vge{\varepsilon}
\def\gf{\phi}       \def\vgf{\varphi}    \def\gh{\eta}
\def\gi{\iota}      \def\gk{\kappa}      \def\gl{\lambda}
\def\gm{\mu}        \def\gn{\nu}         \def\gp{\pi}
\def\vgp{\varpi}    \def\gr{\rho}        \def\vgr{\varrho}
\def\gs{\sigma}     \def\vgs{\varsigma}  \def\gt{\tau}
\def\gu{\upsilon}   \def\gv{\vartheta}   \def\gw{\omega}
\def\gx{\xi}        \def\gy{\psi}        \def\gz{\zeta}
\def\Gg{\Gamma}     \def\Gd{\Delta}      \def\Gf{\Phi}
\def\Gth{\Theta}
\def\Gl{\Lambda}    \def\Gs{\Sigma}      \def\Gp{\Pi}
\def\Gw{\Omega}     \def\Gx{\Xi}         \def\Gy{\Psi}

\def\CS{{\mathcal S}}   \def\CM{{\mathcal M}}   \def\CN{{\mathcal N}}
\def\CR{{\mathcal R}}   \def\CO{{\mathcal O}}   \def\CP{{\mathcal P}}
\def\CA{{\mathcal A}}   \def\CB{{\mathcal B}}   \def\CC{{\mathcal C}}
\def\CD{{\mathcal D}}   \def\CE{{\mathcal E}}   \def\CF{{\mathcal F}}
\def\CG{{\mathcal G}}   \def\CH{{\mathcal H}}   \def\CI{{\mathcal I}}
\def\CJ{{\mathcal J}}   \def\CK{{\mathcal K}}   \def\CL{{\mathcal L}}
\def\CT{{\mathcal T}}   \def\CU{{\mathcal U}}   \def\CV{{\mathcal V}}
\def\CZ{{\mathcal Z}}   \def\CX{{\mathcal X}}   \def\CY{{\mathcal Y}}
\def\CW{{\mathcal W}} \def\CQ{{\mathcal Q}} 
\def\BBA {\mathbb A}   \def\BBb {\mathbb B}    \def\BBC {\mathbb C}
\def\BBD {\mathbb D}   \def\BBE {\mathbb E}    \def\BBF {\mathbb F}
\def\BBG {\mathbb G}   \def\BBH {\mathbb H}    \def\BBI {\mathbb I}
\def\BBJ {\mathbb J}   \def\BBK {\mathbb K}    \def\BBL {\mathbb L}
\def\BBM {\mathbb M}   \def\BBN {\mathbb N}    \def\BBO {\mathbb O}
\def\BBP {\mathbb P}   \def\BBR {\mathbb R}    \def\BBS {\mathbb S}
\def\BBT {\mathbb T}   \def\BBU {\mathbb U}    \def\BBV {\mathbb V}
\def\BBW {\mathbb W}   \def\BBX {\mathbb X}    \def\BBY {\mathbb Y}
\def\BBZ {\mathbb Z}

\def\GTA {\mathfrak A}   \def\GTB {\mathfrak B}    \def\GTC {\mathfrak C}
\def\GTD {\mathfrak D}   \def\GTE {\mathfrak E}    \def\GTF {\mathfrak F}
\def\GTG {\mathfrak G}   \def\GTH {\mathfrak H}    \def\GTI {\mathfrak I}
\def\GTJ {\mathfrak J}   \def\GTK {\mathfrak K}    \def\GTL {\mathfrak L}
\def\GTM {\mathfrak M}   \def\GTN {\mathfrak N}    \def\GTO {\mathfrak O}
\def\GTP {\mathfrak P}   \def\GTR {\mathfrak R}    \def\GTS {\mathfrak S}
\def\GTT {\mathfrak T}   \def\GTU {\mathfrak U}    \def\GTV {\mathfrak V}
\def\GTW {\mathfrak W}   \def\GTX {\mathfrak X}    \def\GTY {\mathfrak Y}
\def\GTZ {\mathfrak Z}   \def\GTQ {\mathfrak Q}

\font\Sym= msam10 
\def\SYM#1{\hbox{\Sym #1}}
\newcommand{\bdw}{\prt\Gw\xspace}
\medskip
\mysection {Introduction}

\setcounter{equation}{0}

The study of isoslated singularities of solutions of quasilinear equations started with the celebrated works of Serrin \cite{Se1}\cite{Se2} dealing with expressions such as
\begin{equation}\label{gen0}
{\rm div} A\left((x,u,Du)\right)+B(x,u,Du)=0
\end {equation}
where $A$ and $B$ are respectively vector valued and real valued Caratheodory functions satisfying the same power $p$-growth with $p\geq 1$. One of the main results of these works stated that the type of singularities is dictated by the diffusion operator $A$. Later on the particular cases of superlinear semilinear elliptic equations was considered, either with an absorption
 \begin{equation}\label{gen1}
-\Gd u+\abs u^{q-1}u=0
\end {equation}
\cite{BrV}, \cite{Ve0}, or with a source reaction
 \begin{equation}\label{gen2}
\Gd u+u^q=0
\end {equation}
\cite{Li}, \cite{GS}, \cite{Av}, and in all cases $q>1$. One of the main facts of these studies relied in the existence of critical thresholds where the interaction of the diffusion and the reaction terms could create unexpected phenomena. As a natural generalisation, the same analysis was carried on for 
 \begin{equation}\label{gen3}
 -{\rm div}\left(\abs {Du}^{p-2}Du\right)+\abs u^{q-1}u=0
\end {equation}
\cite{FV}, and
 \begin{equation}\label{gen4}
{\rm div}\left(\abs {Du}^{p-2}Du\right)+u^q=0
\end {equation}
\cite {SZ}, in the range $0<p-1<q$. In all these works, the   radial explicit solutions, whenever they exist, played a key role. \smallskip

Similarly, the study of the boundary behaviour of solutions of quasilinear equations, has a natural starting point in the description of their isolated singularities on the boundary. Besides the historical results of Fatou, Herglotz and Doob on the boundary trace of posi
tive harmonic and super harmonic functions, equations of types (\ref{gen1}), (\ref{gen2}) and (\ref{gen3}) have alredy been considered (\cite{BVPV},\cite{GV},\cite{BoV}). In the present article we consider equations of type (\ref{gen4}). The problem can be stated under the following form:
Assume  $\Gw$ is an open subset of $\BBR^N$, $a\in\prt\Gw$ and $u\in C(\overline\Gw\setminus\{a\})\cap C^1(\Gw)$ is a solution of one of the above equations which vanishes on $\prt\Gw\setminus\{a\}$, what is the behaviour of $u(x)$ when $x\to a$. The simplest configuration corresponds to $\Gw=\BBR^N_{+}$, and $a=0$ (or more generaly, if $\Gw$ is a cone and the singular point $a$ its vertex $0$). For such geometry, the key-stone element for describing the behaviour of $u$ near $0$ is played by {\it separable solutions}, whenever they exist.
These solutions, which have the form
\begin{equation}\label{mainI-2}
u(x)=u(r,\gs)=r^{-\gb}\gw(\gs)\quad r>0,\,\gs\in S^{N-1},
\end {equation}
have already proved their importance for (\ref{gen1}), (\ref{gen2}) and (\ref{gen3}). It is expected that such will be the case for (\ref{gen4}), even if the full theory will be much more difficult to develop because of the absence of comparison principle and a priori estimates near $x=0$. It is straightforward that, if $u$ is a separable solution of (\ref{gen4}) in $\BBR^N$, 
\begin{equation}\label{mainI-3}
\gb=\myfrac{p}{q+1-p}:=\gb_q,
\end {equation}
which is positive since $q>p-1$. Furthermore
$\gw$ is a solution of 
\begin{equation}\label{mainI-4}\BA {l}
-\nabla'.\left(\left(\gb_q^2\gw^2+\abs{\nabla'\gw}^2\right)^{p/2-1}\nabla'\gw\right)-
\abs{\gw}^{q-1}\gw=\gl_{q,p}\left(\gb_q^2\gw^2+\abs{\nabla'\gw}^2\right)^{p/2-1}\gw,
\EA\end {equation}
in $S_{+}^{N-1}$, where $\nabla'$ is the covariant gradient on $S^{N-1}$, $\nabla'{\bf.}$ the divergence operator acting on vector fields on $S^{N-1}$ and
$$\gl_{q,p}=\gb_q(q\gb_q-N).
$$ 
When $p=2$, $\gb_{q}=2/(q-1)$ and (\ref{mainI-4}) becomes
\begin{equation}\label{p=2-1}
-\Gd '\gw-\abs\gw^{q-1}\gw=\gl_{q,2}\gw,
\end {equation}
where $\Gd '$ is the Laplace-Beltrami operator on $S^{N-1}$ and
$$\gl_{q,2}=\myfrac{2}{q-1}\left(\myfrac{2q}{q-1}-N\right).
$$
If  $S$ is a subdomain of $ S^{N-1}$, equation (\ref{p=2-1}), considered in $S$, is the Euler-Lagrange variation of the functional
\begin{equation}\label{p=2-J}
 I(\psi)=\myint{S}{}
 \left(\myfrac{1}{2}\abs{\nabla\psi}^2+\myfrac{\gl_{q,2}}{2}\psi^2
-\myfrac{1}{q+1}\abs\psi^{q+1}\right)d\gs.
\end {equation}
For any $1<q<(N+1)/(N-3)$ (any $q>1$ if $N=2$ or $3$) this functional satisfies the Palais-Smale condition. Furthermore, if $\gl_{q,2}<\gl_{S,2}$, ($\gl_{S,2}$ is the first eigenvalue of $-\Gd '$ in $W^{1,2}_{0}(S)$), Ambrosetti-Rabinowitz theorem \cite {AR} or Pohozaev fibration method \cite {Po1}, \cite{Po2} apply and yield to the existence of non-trivial positive solutions to (\ref{p=2-1}) in $S$ vanishing on $\prt S$; while if $\gl_{q,2}\geq\gl_{S,2}$ no such solution exists.\\

When $p\neq 2$, equation (\ref{mainI-4}) cannot be associated to any functional defined on 
$S^{N-1}$, except if $q=q_{c}=(N(p-1)+p)/(N-p)$ (the critical Sobolev exponent for $W^{1,p}$, when $N>p$); therefore, finding functions satisfying it is not straightforward. Besides the constant solutions 
 which exist as soon as $q\gb_q<N$, it is not easy to prove the existence of non-constant solutions.  As in the case $p=2$, it is remarkable to see that existence, or non existence, of solutions of  (\ref{mainI-4}) is associated to some spectral problem, although this problem is not standard at all: if one looks for the existence of a positive p-harmonic function $v$ in the cone 
$C_S=\{(r,\gs):r>0,\gs\in S\}$ vanishing on $\prt S$, under the form $v(r,\gs)=r^{-\gb}\phi(\gs)$, one finds that $\phi$  is a positive solution of the so-called {\it spherical $p$-harmonic spectral equation on $S$}, namely
 \begin{equation}\label{mainI-5}\left\{\BA {l}
-\nabla'.\left(\left(\gb^2\phi^2+\abs{\nabla'\phi}^2\right)^{p/2-1}\nabla'\phi\right)
=\gl\left(\gb^2\phi^2+\abs{\nabla'\phi}^2\right)^{p/2-1}\phi\quad \rm in\; S\\
\phantom{,,,,,,,,+\gl\left(\gb^2\phi^2+\abs{\nabla'\phi}^2\right)^{p/2-1}}
\phi=0\quad \rm in\; \partial S,
\EA\right.\end {equation}
and
$\gl=\gb\left(\gb (p-1)+p-N\right)$. The difficulty of this problem is two-fold since $\gb$ is unknown and (\ref{mainI-5}) is not the Euler-Lagrange equation of any functional. However, 
 given a smooth subdomain $S\subset S^{N-1}$, it is proved in \cite{Ve1}, following a shooting method due to Tolksdorff \cite {To}, that there exists a  couple $(\gb,\phi)=(\gb_S,\phi_S)$, where $\gb_S>0$ is unique and 
 $\phi_S$ is defined up to an homothethy, such that (\ref{mainI-5}) holds. Denoting
 $$\gl_{S}=\gb_{S}\left(\gb_{S} (p-1)+p-N\right),
 $$
the couple $(\phi_S, \gl_{S})$ is the natural generalization of the first eigenfunction and eigenvalue of the Laplace-Beltrami operator in $W^{1,2}_{0}(S)$ since $\gl_{S}=\gl_{S,2}$ when $p=2$. Our first theorem is a non-existence which extends the one already mentioned in the case $p=2$.\\
 
 \noindent {\bf Theorem 1.}{ \it Let $S\subset S^{N-1}$ be a   smooth subdomain. If  $\gb_q\geq\gb_S$ there exists no positive solution of (\ref{mainI-4})
 in $S$ which vanishes on $\prt S$.}\\
 
 Apart the case $p=2$, the existence counterpart of this theorem is not known in arbitrary dimension, except if 
 $q=q_{c}$ in which case (\ref{gen4}) is the Euler-Lagrange equation of the functional
 \begin{equation}\label{q=2-J}
 J(\psi)=\myint{S}{}
 \left(\myfrac{1}{p}\left(\gb^2_{q_{c}}\psi^2+\abs{\nabla'\psi}^2\right)^{p/2}
 -\myfrac{1}{q_{c}+1}\abs\psi^{q_{c}+1}\right)d\gs,
\end {equation}
and applications of the already mentioned variational methods lead to an existence result.\medskip 
 
However, when  $N=2$ the problem of finding solutions of (\ref{gen4}) under the form (\ref{mainI-2}) can be completely solved using dynamical systems methods. In order to point out a richer class of phenomena, we shall imbed this problem into a more general class of quasilinear equations with a potential, authorizing even the value $p=1$. This equation is the following,
 \begin{equation}\label{2dim1}
 {\rm div}\left(\abs {Du}^{p-2}Du\right)+\abs u^{q-1}u-\myfrac{c}{\abs x^p}\abs u^{p-2}u=0
 \end {equation}
 in $\BBR^2\setminus\{0\}$, with $q>p-1\geq 0$ and $c\in\BBR$. If $u$ is a solution under the form (\ref{mainI-2}), $\gb$ is be equal to $\gb_q$, while $\gw$ is any $2\gp$-periodic solution of 
 \begin{equation}\label{2dim2}\BA {l}
\myfrac{d}{d\sigma}\left[  \left(  \beta_{q}^{2}\omega^{2}+\left(\myfrac{d\omega
}{d\sigma}\right)^{2}\right)  ^{p-2)/2}\myfrac{d\omega}{d\sigma}\right]
+\lambda_{q}\,\left[  \beta_{q}^{2}\omega^{2}+\left(\myfrac{d\omega}{d\sigma}%
\right)^{2}\right]  ^{(p-2)/2}\omega\\[4mm]
\quad\qquad\qquad\qquad\qquad\qquad\qquad\qquad\qquad\qquad\qquad\qquad+|\omega|^{q-1}\omega-c|\omega
|^{p-2}\omega=0,%
\EA\end{equation}
where
 \begin{equation}\label{2dim3}
\lambda_{q}=\beta_{q}\,(q\beta_{q}-2)=\beta_{q}\,(p-2+(p-1)\beta_{q}).
\end{equation}
If we set
 \begin{equation}\label{2dim4}
c_{q}=\beta_{q}^{p-2}\lambda_{q}=p^{p-1}\frac{(p-2)q+2(p-1)}{\left(
q+1-p\right)  ^{p}},
\end{equation}
then, if $c\leq c_q$, the only constant solution is the zero function, while
if $c>c_{q}$, there exist two other constant
solutions $\pm(c-c_{q})^{1/(q+1-p)}$. Let us denote by  $\mathcal{E}%
^{+}$ the set of positive solutions of (\ref{2dim2}) on $S^{1}$, $\mathcal{E}$ the set of sign changing solutions and $\CF=\pm\CE^+\cup\CE$ the set of all nonzero solutions.
Our main result which gives the struture of the sets $\CE$ and $\CE^+$ is the following:\\

\noindent{\bf Theorem 2.} {\it Assume $p>1,$ $q>p-1.$ Then\smallskip

\noindent(i)
 \begin{equation}\label{2dim5}
\CE={\displaystyle\bigcup\limits_{\substack{k\in\mathbb{N},\\k=k_{q}}}^{\infty}}
\left\{  \omega_{k}(.+\psi):\psi\in S^{1}\right\},%
\end{equation}
in which expression $\omega_{k}$ is a function with least period $2\pi/k$, and $k_{q}=1$
if $c\geq c_{q}$, or  $k_{q}$ is the smallest positive integer
such that $k_{q}>M_{q},$ where
 \begin{equation}\label{2dim6}
M_{q}=\frac{\pi\beta_{q}^{1-p}}{2%
{\displaystyle\int\nolimits_{0}^{\pi/2}}
\myfrac{1+(p-1)\tan^{2}\theta}{\beta_{q}^{p}(p-1)\tan^{2}\theta+c_{q}%
-c\cos^{p-2}\theta}d\theta},
\end{equation}\smallskip
if $c<c_{q}$.\smallskip

\noindent(ii) If $c\leq c_{q}$, $\mathcal{E}^{+}$ is empty. If $0<c-c_{q}\leq$
$\beta_{q}^{p-1}/p,$ $\mathcal{E}^{+}$ is reduced to the constant function $(c-c_{q})^{1/(q+1-p)}$. If $c-c_{q}>$ $\beta_{q}^{p-1}/p,$ $\mathcal{E}^{+}$ contains the constant function $(c-c_{q})^{1/(q+1-p)}$
and the set 
 \begin{equation}\label{2dim7}
\CE^+_{*}=\bigcup\limits_{\substack{k\in\mathbb{N},\\k=1}}^{k_{q}^{+}}
\left\{  \omega_{k}^{+}(.+\psi):\psi\in S^{1}\right\},
\end{equation}
where $\omega_{k}^{+}$ is a non-constant positive function with least period $2\pi/k,$ and
$k_{q}^{+}$ is the largest integer smaller than $(p\beta_{q}^{1-p}%
(c-c_{q}))^{1/2}.$}\\

Since separable solutions of (\ref{gen4}) defined in a cone $C_{S}$ and vanishing on $\prt C_{S}$ are associated to elements of $\CE$, we can prove the existence counterpart of Theorem 1 in dimension 2.\\

\noindent{\bf Corollary 1. }{\it  Let $N=2$ and $S$ be any angular sector of $S^1$. Then there exists a positive solution of (\ref{mainI-4}) vanishing at the two end points of $S$ if and only if $\gb_{q}<\gb_{S}$. Furthermore this solution is unique. In particular, existence holds for any sector if $p<2$ and $q\geq 2(p-1)/(2-p)$.}\\

The case $p=1$ appears as a limiting case of the preceding one. In that case we observe that $u$ is a positive solution of (\ref{2dim1}) if and only if $v=u^q$ is a solution of the same equation relative to $q=1$,
 \begin{equation}\label{2dim1'}
{\rm div}\left(\abs {Dv}^{-1}Dv\right)+v-\myfrac{c}{\abs x}=0.
 \end {equation}
 The initial case $c=0$ is easily treated, but the case $c\neq 0$, that we shall analyse in full generality, is much richer and delicate and shows a large variety of solutions depending on various parameters.
 \\

\noindent{\bf Theorem 3. }{\it  Assume $p=1$ and $q>0$. Then\smallskip

\noindent(i) If $c\neq0,$ or $c=0$ and $q>1,$ $\mathcal{E}$ is empty. If
$c=0$ and $q\leq1$,
$
\mathcal{E}=\left\{  \omega_{0}(.+\psi):\psi\in S^{1}\right\},
$
where $\gs\mapsto \omega_{0}(\sigma):=2^{1/q}\left\vert \sin\sigma\right\vert ^{(1-q)/q}%
\sin\sigma$ is a $C^{1}$ solution  of (\ref{2dim2}).\smallskip

\noindent(ii) If $c\leq-1$, $\mathcal{E}^{+}$ is empty. If $-1<c<0,$ 
$\mathcal{E}^{+}$ is reduced to the constant function $(c+1)^{1/q}$.  If $c>0$, 
\[
\mathcal{E}^{+}=\left\{  (c+1)^{1/q}\right\}  \cup%
{\displaystyle\bigcup\limits_{\substack{k\in\mathbb{N},\\k=k_{1}}}^{k_{2}}}
\left\{  \omega_{k}^{+}(.+\psi):\psi\in S^{1}\right\} ,
\]
in which expression $\omega_{k}^{+}$ is a positive function with least period $2\pi/k$,
$k_{2}$ is the largest integer strictly smaller than $(c+1)^{1/2}$ and
$k_{1}$ is the smallest integer greater than $\myfrac{\gp}{2}%
{\displaystyle\int\nolimits_{0}^{\pi/2}}
\sqrt{\frac{\cos\theta}{\cos\theta+2c}}d\theta$.
Finally, if $c=0$,
\[
\mathcal{E}^{+}=\left\{  1\right\}  \cup%
{\displaystyle\bigcup\limits_{K\in\left(  0,1\right)  }}
\left\{  \omega_{K}^{+}(.+\psi):\psi\in S^{1}\right\}  \cup\left\{
\begin{array}
[c]{c}%
\varnothing\qquad\qquad\qquad\qquad\;\;\;\text{ if }q\geq1\\[2mm]
\left\{  \omega_{0}^{+}(.+\psi):\psi\in S^{1}\right\}  \quad\text{if }q<1
\end{array}
\right.
\]
where the functions $\omega_{K}^{+}$ and $\omega_{0}^{+}$ are explicitely given by
\[
\omega_{K}^{+}=\left(  \sqrt{1-K^{2}\sin^{2}\sigma}-K\cos\sigma\right)
^{1/q},\;\;\text{and }\;\omega_{0}^{+}=(2\left\vert \sin\sigma\right\vert )^{1/q}%
\forevery\sigma\in S^{1}.%
\]\\
}

A striking phenomenon is the existence of a 2-parameter family of solutions when $c=0$.\smallskip 
\\

\noindent Our paper is organized as follows: 1- Introduction. 2- The N-dimensional case. 3- The 2-dim dynamical system. 4- The case $p>1$. 5- The case $p=1$. 
\section {The N-dimensional case}
\setcounter{equation}{0}
\subsection{The spherical p-harmonic spectral problem}
If  $p\geq 1$, $\gb>0$ and $\gl\in\BBR$ we denote by $\frak T_{\gb,\gl}$ the operator defined on $C^1(S^{N-1})$ by
\begin {equation}\label {T}
\varphi\mapsto\frak T_{\gb,\gl}[\varphi]=-\nabla'.\left( (\beta^2\varphi^2+{\abs 
{\nabla'\varphi}}^2)^{(p-2)/2}\nabla'\varphi\right)
-\gl (\beta^2\varphi^2+{\abs {\nabla'\varphi}}^2)^{(p-2)/2}\varphi
.
\end {equation}
Let $q>p-1>0$, $S$ be a smooth connected domain on $S^{N-1}$ and $C_S$ the cone with vertex $0$ generated by $S$. If $u$ is 
a positive solutions of  
\begin {eqnarray}\label {singular}
-{\rm div}({\abs{D u}^{p-2}}D u)=u^q,
\end {eqnarray}
in $C_S\setminus\{(0)\}$ vanishing on $\prt C_S\setminus\{(0)\}$, under the form
\begin {eqnarray}\label {pseudorad}
u(r,\gs)=r^{-\gb}\gw (\gs),
\end {eqnarray}
 then 
$\beta=p/(q+1-p):=\beta_{q}$ and $\omega$ solves
\begin  {equation}\label {anisotropic}\left\{  \BA {l}
\frak T_{\gb_{q},\gl_{q,p}}[\gw]-\gw^q
=0
\mbox { in }\,  S\\[2mm]
\phantom {\frak T_{\gb_{q},\gl_{q,p}}[\gw]-\,^q}
\gw=0\, \mbox { on }\, \prt S,
\EA\right.\end {equation}
 where
$$
\gl_{q,p}=\beta_{q}(q\beta_{q}-N).
$$
We denote by $\gb_S$ the exponent corresponding to the first spherical singular $p$-harmonic function and by 
$\phi_S$ the corresponding function. Thus $\gb_S>0$ and $u(r,\gs)=r^{-\gb_S}\phi_S (\gs)$ is $p$-harmonic
in $C_S\setminus\{(0)\}$ and vanishes on $\prt C_S\setminus\{(0)\}$. Furthermore $\phi=\phi_S>0$ and satisfies
\begin  {equation}\label {anisotropic2}\left\{  \BA {l}
\frak T_{\gb_{S},\gl_{S}}[\phi]=0\;
\mbox { in }\;  S\\[2mm]
\phantom {\frak T_{\gb_{S},\gl_{S}}[]}
\phi=0\, \mbox { on }\, \prt S,
\EA\right.
\end {equation}
 where
$$\gl_S=\gb_S(\gb_S(p-1)+p-N).
$$

We recall that $(\gb_S,\phi_S)$ is unique up to an homothety upon $\phi$. Furthermore 
$\phi_S$ is positive in $S$, $\prt\phi_S/\prt\gn<0$ on $\prt S$ and
$$S'\subset S,\;S'\neq S\Longrightarrow \gb_{S'}>\gb_S.$$
\subsection{Non-existence}

\noindent{\it Proof of Theorem 1.} We  put 
$$\gth=\myfrac {\gb_q}{\gb_{S}}\;\mbox { and }\;\eta=\phi_{S}^\gth.
$$
Then $\gth\geq 1$ and
$$\nabla' \eta=\gth\phi_{S}^{\gth-1}\nabla' \phi_{S},
$$
$$\gb^2_q\eta^2+\abs {\nabla' \eta}^2=\gth^{2}\phi_{S}^{2(\gth-1)}
(\gb^{2}_{S}\phi_{S}^{2}+\abs {\nabla' \phi_{S}}^2),
$$
$$(\gb^2_q\eta^2+\abs {\nabla' \eta}^2)^{(p-2)/2}=
\gth^{p-2}\phi_{S}^{(p-2)(\gth-1)}(\gb^{2}_{S}\phi_{S}^{2}+\abs {\nabla' \phi_{S}}^2)^{(p-2)/2},
$$
$$\BA {l}\nabla'.(\gb^2_q\eta^2+\abs {\nabla' \eta}^2)^{(p-2)/2}\nabla'\eta
=\gth^{p-1}\phi_{S}^{(p-1)(\gth-1)}
\nabla'.(\gb^{2}_{S}\phi_{S}^{2}+\abs {\nabla' \phi_{S}}^2)^{(p-2)/2}\nabla' \phi_{S}\\[2mm]
\phantom {------}+(p-1)(\gth-1)
\gth^{p-2}\phi_{S}^{(p-1)(\gth-1)-1}(\gb^{2}_{S}\phi_{S}^{2}+\abs {\nabla' \phi_{S}}^2)^{(p-2)/2}\abs {\nabla' \phi_{S}}^2
\EA$$
Using (\ref {anisotropic2}) with $\phi=\phi_{S}$, we derive
\begin  {equation}\label {anisotropic3}\frak T_{\gb_{q},\gl_{q,p}}[\eta]
=-(p-1)\gth^{p-1}(\gth-1)\phi_{S}^{(p-1)(\gth-1)-1}(\gb^{2}_{S}\phi_{S}^{2}+\abs {\nabla' \phi_{S}}^2)^{p/2}
\mbox { in }\,  S.
\end {equation}
Because $\gw$ is a nonnegative nontrivial solution of (\ref {anisotropic}), it is nonpositive in $S$. Furthermore $\prt\gw/\prt\gn<0$ on $\prt S$. Therefore we can choose $\phi_S$ as the maximal positive solution of (\ref {anisotropic2}) such that $\eta\leq \gw$. If $\gth>1$ there exists $\gs^*\in S$ such that 
\begin  {equation}\label {tangent1}
\gw(\gs^*)=\eta(\gs^*)>0\;\mbox { and }\;\gw(\gs)\geq\eta(\gs)\forevery
\gs\in \bar S.
\end {equation}
If $\gth=1$, the graphs of $\gw$ and $\eta$ could be tangent only on $\prt S$. This means that either (\ref {tangent1}) holds, or there exists $\bar\gs\in\prt S$ such that
\begin  {equation}\label {tangent2}
\prt\gw(\bar\gs)/\prt\gn=\prt\eta(\bar\gs)/\prt\gn<0\;\mbox { and }\;\gw(\gs)<\eta(\gs)\forevery
\gs\in S.
\end {equation}
Let $\psi=\gw-\eta$ and we first consider the case where (\ref {tangent1}) holds. Let $g=(g_{ij})$ be the 
metric tensor on $S^{N-1}$. We recall the following expressions  in local 
coordinates $\sigma_{j}$ around $\gs^*$, 
$${\abs {\nabla'\varphi }^{2}}=
\sum_{j,k}g^{jk}\frac {\partial \varphi}{\partial \sigma_{j}}\frac {\partial 
\varphi}{\partial \sigma_{k}},
$$
for any $\varphi\in C^1(S)$, and 
$$\nabla'.X=\frac {1}{\sqrt {\abs g}} \sum_{\ell}
\frac {\partial}{\partial \sigma_{\ell}}\left(\sqrt {\abs g}X^\ell\right)
=\frac {1}{\sqrt {\abs g}} \sum_{\ell,i}\frac {\partial}{\partial \sigma_{\ell}}\left(\sqrt {\abs g}g^{\ell i}X_{i}\right)
,$$
for any vector field $X\in C^1(TS^{N-1})$, if we lower the indices by setting 
$\displaystyle {X^\ell=\sum_{i}g^{\ell i}X_{i}}$. We derive from the mean value theorem
\begin {eqnarray*}(\beta^{2}_{q}\omega^2+\abs 
{\nabla'\omega}^2)^{(p-2)/2}\frac {\partial \omega}{\partial 
\sigma_{i}}-
(\beta^{2}_{q}\eta^2+\abs 
{\nabla'\eta}^2)^{(p-2)/2}\frac {\partial \eta}{\partial 
\sigma_{i}}
=\sum_{j}\alpha^i_{j}\frac {\partial (\omega-\eta)}{\partial 
\sigma_{j}}+b^i(\omega-\eta),
\end {eqnarray*}
where 
\begin {eqnarray*}b^i=(p-2)\left(\beta_{q}^2(\eta+t(\omega-\eta))^2+{\abs 
{\nabla'(\eta+t(\omega-\eta))}}^2\right)^{(p-4)/2}\\
\times (\eta+t(\omega-\eta))\frac {\partial 
(\eta+t(\omega-\eta))}{\partial\sigma_{i}},
\end {eqnarray*}
and
\begin {eqnarray*}
\alpha^i_{j}=(p-2)\left(\beta_{q}^2(\eta+t(\omega-\eta))^2+{\abs 
{\nabla'(\eta+t(\omega-\eta))}}^2\right)^{(p-4)/2}\\
\times\frac {\partial 
(\eta+t(\omega-\eta))}{\partial\sigma_{i}}\sum_{k}g^{jk}\frac {\partial 
(\eta+t(\omega-\eta))}{\partial\sigma_{k}}\\
+\delta_{i}^j\left(\beta_{q}^2(\eta+t(\omega-\eta))^2+{\abs 
{\nabla'(\eta+t(\omega-\eta))}}^2\right)^{(p-2)/2}.
\end {eqnarray*}
Since the graph of $\eta$ and $\omega $ are tangent at 
$\sigma^*$, 
$$\eta(\sigma^*)=\omega(\sigma^*)=P_{0}>0\quad \mbox {and } 
\nabla'{\eta(\sigma^*)}=\nabla'{\omega(\sigma^*)=Q}.$$
Thus
$$b^i(\sigma^*)=
(p-2)\left(\beta_{q}^2P_{0}^2+{\abs 
{Q}}^2\right)^{(p-4)/2}P_{0}Q_{i},
$$
and
\begin {eqnarray*}
\alpha^i_{j}(\sigma^*)=\left(\beta_{q}^2P_{0}^2+{\abs Q}^2\right)^{(p-4)/2}
\left(\delta_{i}^j(\beta_{q}^2P_{0}^2+{\abs Q}^2)+(p-2) Q_{i}\sum_{k}g^{jk}Q_{k}\right).
\end {eqnarray*}
Now
\begin {eqnarray*}
\frak T_{\gb_{q},\gl_{q,p}}[\gw]-\frak T_{\gb_{q},\gl_{q,p}}[\eta]=\gw^q+(p-1)\gth^{p-1}(\gth-1)\phi_{S}^{(p-1)(\gth-1)-1}(\gb^{2}_{S}\phi_{S}^{2}+\abs {\nabla' \phi_{S}}^2)^{p/2}
\phantom {------}
\\
 =\frac {-1}{\sqrt {\abs g}}\sum_{\ell,i}\frac {\partial}{\partial \sigma_{\ell}}\left[\sqrt {\abs g}g^{\ell i}
\left((\beta^{2}_{q}\gw^2+\abs 
{\nabla'\gw}^2)^{\frac {p}{2}-1}\frac {\partial \gw}{\partial 
\sigma_{i}}-
(\beta^{2}_{q}\eta^2+\abs 
{\nabla'\eta}^2)^{\frac {p}{2}-1}\frac {\partial \eta}{\partial 
\sigma_{i}}\right)\right]\\
-\gl_{q,p}\left((\beta^{2}_{q}\gw^2+\abs 
{\nabla'\gw}^2)^{\frac {p}{2}-1}\gw-(\beta^{2}_{q}\eta^2+\abs 
{\nabla'\eta}^2)^{\frac 
{p}{2}-1}\eta\right),\\
=-\frac {1}{\sqrt {\abs g}}\sum_{\ell,i}\frac {\partial}{\partial \sigma_{\ell}}\left[\sqrt {\abs g}g^{\ell i}
\left(\sum_{j}\alpha^i_{j}\frac {\partial (\gw-\eta)}{\partial \sigma_{j}}
+b^i(\gw-\eta)\right)\right]
+\sum_{i}C_{i}\frac {\partial (\gw-\eta)}{\partial 
\sigma_{i}}\\
=-\frac {1}{\sqrt {\abs g}}\sum_{\ell,j}\frac {\partial}{\partial 
\sigma_{\ell}}\left[a^\ell_{j}\frac {\partial (\gw-\eta)}{\partial \sigma_{j}}\right]
+\sum_{i}C_{i}\frac {\partial (\gw-\eta)}{\partial 
\sigma_{i}},\phantom {-----------.-}
\end {eqnarray*}
where the $C_{i}$ are continuous functions and 
$$a^\ell_{j}=\sqrt{\abs g}\sum_{i}g^{\ell i}\alpha^i_{j}.
$$
The matrix $\left(\alpha^i_{j}(\sigma_{0})\right)$ is symmetric, 
definite and positive since it is the Hessian of the strictly convex function
$$X=(X_{1},\ldots,X_{n-1})\mapsto \frac {1}{p}\left(P_{0}^{2}+{\abs X}^2\right)^{p/2}
=\frac {1}{p}\left(P_{0}^{2}+\sum_{j,k}g^{jk}X_{j}X_{k}\right)^{p/2}.
$$
Therefore $\left(\alpha^i_{j}\right)$ has the same property in 
some neighborhood of $\gs^*$, and 
the same holds true with $\left(a^\ell_{j}\right)$. Finally the 
function $\psi=\omega-\eta$ is nonnegative, vanishes at $\gs^*$ 
and satifies
\begin {eqnarray}\label {max princ}-\frac {1}{\sqrt {\abs g}}\sum_{\ell,j}\frac {\partial}{\partial 
\sigma_{\ell}}\left[a^\ell_{j}\frac {\partial \psi}{\partial \sigma_{j}}\right]
+\sum_{i}C_{i}\frac {\partial \psi}{\partial 
\sigma_{i}}\geq 0.
\end {eqnarray}
Then $\psi=0$ in a neighborhood of $S$. Since $S$ is connected, $\psi$ 
is identically $0$ which a contradiction. 

If (\ref {tangent2}) holds, then $\theta=1$ and the graphs of 
$\eta$ and $\omega$ are tangent at  $\bar\sigma$. Proceeding as above and using the fact 
that $\partial \eta/\partial\nu$ exists and never vanishes on the 
boundary, we see that $\psi=\eta-\omega$ satisfies (\ref {max 
princ}) with a strongly elliptic operator in a neighborhood $\mathcal 
N$ of $\bar\sigma$. Moreover $\psi >0$ in $\mathcal N$, $\psi 
(\bar\sigma)=0$ and $\partial \psi/\partial\nu(\bar\sigma)=0$. This 
is a contradiction, which ends the proof.\qeda\\

\noindent \Remark If $p=2$, the proof of non-existence is straightforward by multiplying the equation in $\gw$ by the first eigenfunction $\phi_S$ and get
$$\myint{S}{}\left((\gl_S-\gl_{q,2})\gw-\gw^q\right)\phi_S d\gs=0,
$$
a contradiction since $\gl_S\leq \gl_{q,2}$. 
\subsection{Existence results}

Let us consider the case $q=q_{c}=(N(p-1)+p)/(N-p)$ ($N>p>1$), and let $S$ be any smooth subdomain of $S^{N-1}$. Since in that case $\gl_{q,p}=-\gb^2_{q_{c}}$, the research of solutions of (\ref{gen4}) under the form (\ref{mainI-2}) vanishing on $\prt C_{S}$ leads to 
\begin{equation}\label {q=3-J}\left\{\BA {l}
\frak T_{\gb_{q_{c}},-\gb^2_{q_{c}}}[\gw]-
\abs{\gw}^{q_{c}-1}\gw=0\quad \mbox {in } S\\[2mm]
\phantom{\frak T_{\gb_{q_{c}},-\gb^2_{q_{c}}}[\gw]-
\abs{\gw}^{q_{c}-1}}
\gw=0\quad \mbox {in }\prt S,
\EA\right.
\end {equation}
where $\gb_{q_{c}}=N/p-1$. This equation is the Euler-Lagrange variation of the functional $J$ defined on $W^{1,p}_{0}(S)$ by
 \begin{equation}\label{q=2-J'}
 J(\psi)=\myint{S}{}
 \left(\myfrac{1}{p}\left(\gb^2_{q_{c}}\psi^2+\abs{\nabla'\psi}^2\right)^{p/2}
 -\myfrac{1}{q_{c}+1}\abs\psi^{q_{c}+1}\right)d\gs.
\end {equation}
\bth {Exist}Problem (\ref{q=3-J}) admits a positive solution.
\es
\Proof Clearly the functional is well defined on $W^{1,p}_{0}(S)$ since $q_{c}$ is smaller than the Sobolev exponent $p^*_{{N-1}}$ for $W^{1,p}$ in dimension N-1. For any $\psi\in W^{1,p}_{0}(S)$, $\lim_{t\to\infty}J(t\psi)=-\infty$. Furthermore there exist $\gd>0$ and $\ge>0$ such that $J(\psi)\geq\ge$ for any $\psi\in W^{1,p}_{0}(S)$ such that $\norm{\psi}_{W^{1,p}}=\gd$. Assume now that $\{\psi_{n}\}$ is a sequence of $W^{1,p}_{0}(S)$ such that 
$J(\psi_{n})\to \ga$ and $\norm{DJ(\psi_{n})}_{W^{-1,p'}}\to 0$ as $n\to\infty$.
Then
$$\frak T_{\gb_{q_{c}},-\gb^2_{q_{c}}}[\psi_{n}]-
\abs{\psi_{n}}^{q_{c}-1}\psi_{n}=\ge_{n}\to 0.
$$
Then
$$\BA {l}
\myint{S}{}\left(\left(\gb_{q_{c}}^2\psi_{n}^2+\abs{\nabla'\psi_{n}}^2\right)^{p/2}
-\abs{\psi_{n}}^{q_{c}+1}\right)d\gs=\langle\ge_{n},\psi_{n}\rangle.
\EA$$
Since $J(\psi_{n})\to \ga$ it follows 
$$\myint{S}{}\left(\gb_{q_{c}}^2\psi_{n}^2+\abs{\nabla'\psi_{n}}^2\right)^{p/2}d\gs
\to p(q_{c}+1)\ga/(q_{c}+1-p).
$$
Therefore $\{\psi_{n}\}$ remains bounded in $L^{q_{c}+1}(S)$, and relatively compact in $L^r(S)$, for any $1<r<q_{c}+1$. Multiplying the equation $DJ(\psi_{n})-\ge_{n}$ by 
$T_{k,\theta}(\psi_{n})$ where 
$\theta\in (1,(p^*_{{N-1}}-1)/q_{c})$, $k>0$ and $T_{k,\theta}(r)={\rm sgn}\min\{\abs r,k\}$ and using standard bootstrap arguments yields to the boundedness of $\{\psi_{n}\}$ in $L^\infty(S)$. Combining this fact with the compactness of  in $L^r(S)$, we derive the compactness in any any $L^s$, for $s<\infty$. Therefore $\{\psi_{n}\}$ is relatively compact in $W^{1,p}_{0}(S)$. This means that $J$ satisfies the Palais-Smale condition. \qeda 

\section {The 2-dim dynamical system}
\setcounter{equation}{0}
\subsection {Extension of the data}
Due to possible applications and similarly to what is done in the semilinear case $p=2$ (see \cite{BiBo}, \cite{ChDu}, \cite{ChJa}), we shall consider the existence  problem for $2\pi$-periodic solutions of a more general quasilinear equation than (\ref{2dim2}), 
\begin{equation}
\frac{d}{d\sigma}\left[  \left(  \beta^{2}\omega^{2}+\left(\frac{d\omega}{d\sigma
}\right)^{2}\right)  ^{p/2-1}\frac{d\omega}{d\sigma}\right]  +\lambda
\,\left[  \beta^{2}\omega^{2}+\left(\frac{d\omega}{d\sigma}\right)^{2}\right]
^{p/2-1}\omega+g(\omega)-c|\omega|^{p-2}\omega=0, \label{Eab}%
\end{equation}
where $\lambda,\beta,c$ are real parameters, with $\beta>0,$ and
$g\in C^{0}(\mathbb{R})\cap$ $C^{1}(\mathbb{R}\backslash\left\{  0\right\}
) $ is odd and satisfies
\begin{equation}
\lim_{s\rightarrow0+}g(s)/s^{q}=1,\qquad\lim_{s\rightarrow\infty}%
g(s)/s^{p-1}=\infty,\qquad\frac{d}{ds}(g(s)/\left\vert s\right\vert
^{p-1})>0\text{ on }\left(  0,\infty\right), \label{1f}%
\end{equation}
with $q>p-1\geq0$. In fact we can easily reduce the problem to a simpler form, and
particularly in the case $p=1,$ where the equation has a remarkable
homogeneity property. The next statement is a straightforward computation which transforms the equation satisfied by $\gw$ into two more canonic forms.\medskip

\blemma{lem} Let $\gw$ be a solution of (\ref{Eab}).\smallskip

\noindent (i) Assume $p>1.$ If we set
\begin{equation}
\tau=\beta\sigma\;,\text{ }\;\omega\left(  \sigma\right)  =\beta
^{p/(q+1-p)}w(\tau)\;\text{ and } w^{\prime}=\myfrac{dw}{d\tau},\label{cv}%
\end{equation}
then $w$ satisfies%
\begin{equation}
\frac{d}{d\tau}\left(  \left(  w^{2}+w^{\prime2}\right)  ^{p/2-1%
}w^{\prime}\right)  -b\,\left(  w^{2}+w^{\prime2}\right)  ^{p/2-1%
}w+f(w)-d|w|^{p-2}w=0, \label{E}%
\end{equation}
where
\begin{equation}
b=\frac{-\lambda}{\beta^{2}},\qquad d=\frac{c}{\beta^{p}},\qquad f(s)=\beta
^{-pq/(q+1-p)}g(\beta^{p/(q+1-p)}s). \label{bd}%
\end{equation}
In particular $f$ satisfies the same assumptions (\ref{1f}) as $g.$\smallskip

\noindent(ii) Assume $p>1.$ If on any open interval $I\subset (0,2\gp)$ where $\omega(\sigma)\neq0,$ we set
\begin{equation}
\tau=\beta q\;\sigma,\quad\text{and }\omega\left(  \sigma\right)  =(\beta
q)^{1/q}\left\vert w(\tau)\right\vert ^{1/q-1}w(\tau), \label{cvq}%
\end{equation}
then $w$ satisfies (\ref{E}) on $I$, with%
\begin{equation}
b=-\lambda/\beta^{2}q,\qquad d=c/\beta q,\qquad f_{1}(s)=\beta^{-q}g((\beta
qs)^{1/q}). \label{b1d1}%
\end{equation}
Furthermore $f_{1}$ satisfies the assumptions (\ref{1f}) with $q=1$, i.e.%
\begin{equation}
\lim_{s\rightarrow0+}f_{1}(s)/s=1,\qquad\lim_{s\rightarrow\infty}%
f_{1}(s)=\infty,\qquad f_{1}^{\prime}(s)>0\text{ on }\left(  0,\infty\right). \label{hfun}%
\end{equation}
\es

Due to this result, the changes of variables (\ref{cv}) and (\ref{cvq}) reduce the
problem to the study both of existence of periodic solutions of equation (\ref{E}), and to characterizing the period function of these solutions, in the range $q>p-1$ if $p>1,$ and $q>0$ if $p=1.$

\subsection{Reduction to dynamical systems\bigskip}

\bigskip

We re-write (\ref{E})  as the system,
\begin{equation}\left\{\BA {l}
w^{\prime}=F(w,y)=y\\
y^{\prime}=G(w,y)=\myfrac{bw^{3}+(b+2-p)w\,y^{2}%
-(f(w)-d|w|^{p-2}w)(w^{2}+y^{2})^{2-p/2}}{w^{2}+(p-1)y^{2}},\EA\right.
\label{FG}%
\end{equation}
and we denote by $h$ the odd function defined on $\BBR$ by
\begin{equation}
h(s)=\left\{\BA {l}f(s)/\abs s^{p-2}s\quad \;\mbox{if }s\neq 0\\[2mm]
0\phantom{---;,---}\mbox {if }s= 0.
\EA\right.
 \label{1h}%
\end{equation}
If $b+d\leq0,$ (\ref{FG}) has no non-trivial stationary point, while if  $b+d>0,$ it admits the two
stationary points $\pm P_{0},$ with
$P_{0}=(a,0)$ and $a=h^{-1}(b+d).$
Furthermore $P_{0}$ is a center since the linearized system at $P_{0}$ is given by the
matrix $$%
\begin{pmatrix}
0 & 1\\
-ah^{\prime}(a) & 0
\end{pmatrix}
.$$
System (\ref{FG}) is clearly singular at $(0,0).$ Furthermore it could singular be along the
line $w=0$  if $p=1$, if $q<1$, and if $p<2$ and  $d\neq0$. Actually, for $p>1$ it
is not singular at any points $(0,\sigma)$ with $\sigma\neq0$. This can be checked as follows:  consider the Cauchy problem
\begin{equation}\label{cauch}\left\{\BA {l}
w^{\prime\prime}=G(w,w^{\prime}),\quad t\in(-\gd,\gd)\\
w(0)=0,w^{\prime}(0)=\sigma,
\EA\right.\end{equation}
and let $w$ be any local solution; since near $(0,\sigma),$ $\ G$ is
continuous with respect to $w$ and $C^{1}$ with respect to $y,$  $w$ is
$C^{2};$ because $\gs\neq 0$,  $t$ can be expressed locally in terms of $w$.  Defining $w^{\prime}(t)=p(w),$ then $p$ is $C^{1}$ near $0,$ $p(0)=1$ and satisfies 
$$\myfrac{dp}{dw}=\myfrac{G(w,p)}{p},$$ 
with $J(w,p)=G(w,p)/p$. Clearly is $C^{1}$ with respect to $p$ and continuous with respect
to $w,$ thus one gets local uniqueness of $p.$ and then the local uniqueness
of problem $w^{\prime}(t)=p(w(t)),$ $w(0)=1$, since $p$ is of class $C^{1}.$\\

The phase plane of the system (\ref{FG}) is equivariant under symmetries with respect to
the two axes of coordinates, because $F$ is even with respect to $w$ and odd
with respect to $y,$ and $G$ is odd with respect to $w$ and even with respect
to $y.$ Thus from now we can restrict the study to the first quadrant
\[
\overline{\mathcal{Q}}\backslash\left\{  (0,0)\right\} ,\quad\text{where
}\mathcal{Q}=\left(  0,\infty\right)  \times\left(  0,\infty\right),
\]
where, in particular, $w\geq0.$  Due to the symmetries, in the case $p>1,$, any trajectory which meets the two axes in finite times $\tau,\tau+T$ is a closed orbit of
period $4T.$\\

\noindent\Remark It is 
useful to introduce the slope $\xi=w^{\prime}/w,$ (or a function
of the slope) as a new variable. This was first used for $p>1$ in \cite{KM} for
the homogeneous problem
\[
\frac{d}{d\tau}\left(  \left(  w^{2}+w^{\prime2}\right)  ^{p/2-1}w^{\prime}\right)  -b\,\left(  w^{2}+w^{\prime2}\right)^{p/2-1}w=0.
\]
In that case the function $\xi$ satisfies
\[
\myfrac{d}{d\tau}\left(  \left(  1+\xi^{2}\right)  ^{p/2-1}\xi\right)
=-((p-1)\xi^{2}-b)\left(  1+\xi^{2}\right)  ^{p/2-1},
\]
for $w>0$, and this equation is completely integrable in terms of $u=\left(
1+\xi^{2}\right)  ^{p/2-1}\xi.$\\

By using polar coordinates  in $\mathcal{Q}$
\[
(w,y)=(\rho\cos\theta,\rho\sin\theta),\qquad\rho>0,\;\theta\in(0,\pi/2),
\]
 we transform (\ref{FG}) into
\begin{equation}\label{pol}\left\{\BA {l}
\theta^{\prime}=\myfrac{b-(p-1)\tan^{2}\theta+(d-h(\rho\cos\theta))\cos
^{p-2}\theta}{1+(p-1)\tan^{2}\theta}\\
\rho^{\prime}=\rho(1+\theta
^{\prime})\tan\theta. %
\EA\right.\end{equation}
Equivalently, if we introduce the slope $\xi=\tan\theta\in\left(  0,\infty
\right)  $, and set
\begin{equation}
u=\phi(\xi)=\cos^{1-p}\theta\sin\theta,\qquad\phi(\xi)=(1+\xi^{2}%
)^{(p-2)/2}\xi,\label{polu}%
\end{equation}
then $\phi^{\prime}(\xi)=(1+\xi
^{2})^{(p-4)/2}(1+(p-1)\xi^{2});$ thus $\phi$ is strictly increasing: from $\left(  0,\infty\right)  $
into $\left(  0,\infty\right)  $ when $p>1,$ and from $\left(  0,\infty\right)  $
into $\left(  0,1\right)  $ when $p=1.$ Defining
\begin{equation}
\varphi=\phi^{-1},\text{ \quad and \quad}E(\xi)=\left(  (p-1)\xi^{2}-b\right)
(1+\xi^{2})^{p/2-1}, \label{D}%
\end{equation}
we obtain 
\begin{equation}\label{suw}
\left\{
\begin{array}
[c]{l}%
w^{\prime}=w\varphi(u),\\[2mm]
u^{\prime}=-E(\varphi(u))-h(w)+d.
\end{array}
\right.  %
\end{equation}
This system is still singular on the line $w=0$ if $h$ $\not \in $ $C^{1}\left(  \left[
0,\infty\right)  \right)  $ near $0.$ In the sequel we set
\begin{equation}
\Psi(u)=%
{\displaystyle\int\limits_{0}^{u}}
\varphi(s)ds. \label{psu}%
\end{equation}

Noticing that
\begin{equation}\label{DerivE}
E^{\prime}(\xi)=\left(  p(p-1)\xi^{2}+2(p-1)-(p-2)b\right)  (1+\xi
^{2})^{(p-2)/2}\xi,
\end{equation}
we derive that $E$ is increasing on $\left(  0,\infty\right)  $ when $(p-2)b\leq 2(p-1)$.
When $(p-2)b> 2(p-1)$, $E$ is decreasing on $\left(  0,\eta\right)  $ and then
increasing, where $\eta$ is defined by
\begin{equation}
p(p-1)\eta^{2}=(p-2)b-2(p-1), \label{het}%
\end{equation}
and
\begin{equation}
\min E=E(\eta)=-\frac{2}{p-2}\left(  \frac{(p-2)(b+p-1)}{p(p-1)}\right)
^{p/2}. \label{ehe}%
\end{equation}
In the case of initial problem (\ref{2dim2}), $E$ is increasing.\\

\noindent\Remark If $p>1$, system (\ref{FG}) is singular at
$(0,0).$ If we replace the assumption $\lim_{s\rightarrow0+}f(s)/s^{q}=1,$ by
the stronger one%
\begin{equation}
\lim_{s\rightarrow0+}f^{\prime}(s)/s^{q-1}=q, \label{lfp}%
\end{equation}
 we can transform system (\ref{suw}) in $\left(  0,\infty\right)
\times\mathbb{R}$ in a system of the same type, but without singularity: this is obtained by
performing the substitution
$
v=w^{q+1-p}.%
$
Then
\begin{equation}\left\{\BA {l}
v^{\prime}=(q+1-p)v\varphi(u)\\
u^{\prime}=-E(\varphi(u))-\tilde
{h}(v)+d, \label{suvt}%
\EA\right.\end{equation}
where $v\mapsto\tilde{h}(v)=h(v^{1/(q+1-p})\in C^{1}(\left[  0,1\right).$ In
particular, if $f(w)=|w|^{q-1}w$, we find
\begin{equation}\left\{\BA {l}
v^{\prime}=(q+1-p)v\varphi(u)\\
 u^{\prime}=-E(\varphi(u))-v+d.
\label{suv}%
\EA\right.\end{equation}

\noindent\Remark In the case $f(w)=|w|^{q-1}w,$ we can differentiate the equation relative to
$u^{\prime}$ and obtain  that $u$  satisfies the following equation
\begin{equation}
u^{\prime\prime}=B(\varphi(u))u^{\prime}+(q+1-p)(E(\varphi(u))-d)\varphi(u),
\label{Eu}%
\end{equation}
where $E$ is given above, and
\begin{equation}
B(\xi)=\frac{(p-2)b+q-3(p-1)+(q+1-2p)(p-1)\xi^{2}}{1+(p-1)\xi^{2}}\xi.
\label{bks}%
\end{equation}
Notice that equation (\ref{Eu}) has no singularity for $p>1.$
\section{The case $p>1$}
\setcounter{equation}{0}
\subsection{Existence of a first integral}

A natural question is to see if equation (\ref{E}) admits a variational
structure.
When $p=2$, it is the case, for any $b$ and $d$. Since (\ref{E}) takes the
form 
$$w^{\prime\prime}-(b\,+d)w+f(w)=0,$$ 
it is the Euler equation of the functional
\[
\mathcal{H}_{2}\mathcal{(}w,w^{\prime})=\frac{w^{\prime2}}{2}+(b+d)\frac
{w^{2}}{2}-\mathcal{F}(w)
\]
where  $\mathcal{F}(w)=
{\displaystyle\int\nolimits_{0}^{w}}f(s)ds.$ 
Thus the function $w^{\prime2}=(b\,+d)w^{2}-2\mathcal{F}%
(w)$ is constant along the trajectory.
When \textbf{ }$p\neq2,p>1$, we find that a first integral exists only in the case $b=1.$
In such a case (\ref{E}) is the Euler equation of the functional
\[
\mathcal{H(}w,w^{\prime})=\frac{\left(  w^{2}+w^{\prime2}\right)  ^{p/2}}%
{p}+d\frac{\left\vert w\right\vert ^{p}}{p}-\mathcal{F}(w).
\]
Therefore, the associated Painlev\'{e} integral
\begin{equation}
\CP(w,w')=\frac{1}{p}\left(  w^{2}+w^{\prime2}\right)  ^{p/2-1}\left(
(p-1)w^{\prime2}-w^{2}\right)  -\frac{d\left\vert w\right\vert ^{p}}%
{p}+\mathcal{F}(w) \label{pain}%
\end{equation}
 is constant along the trajectories. Using the function $E$ introduced at (\ref{D}), then (\ref{pain}) is
equivalent to
\begin{equation}
E\left(\frac{w^{\prime}}{w}\right)=E(\varphi(u))=d-p\frac{K+\mathcal{F}(w)}{w^{p}}
\label{emo}%
\end{equation}
 for $w>0$.  Hence $E$ is increasing on $\left(  0,\infty\right)  $ from $-b=-1$ to
$+\infty$.\medskip\ 

In the general case, we cannot use a first integral for studying the
periodicity properties of the solutions, while it was the main tool in \cite{BiBo} for
$p=2$. This is the reason for which we are lead to use phase plane techniques. Notice that, for the initial problem
(\ref{2dim2}), the value $b=1$ corresponds to the case $p<2$ and
$q=(3p-2)/(2-p)).$
\subsection{Description of the solutions}

In this section we describe in full details the trajectories of system (\ref{FG}) in the
phase plane $(w,y)$. Notice  that the system can be singular on the axis $w=0.$ 

\bprop{Pr1}
Assume $p>1.$ Then all the orbits of system (\ref{FG}) are bounded. Any
trajectory $\mathcal{T}_{\left[  P\right]  }$ issued from a point $P$ in
$\mathcal{Q}$ is \smallskip

\noindent  (i) either a closed orbit surrounding $(0,0)$,  \smallskip

\noindent (ii) or, if
$b+d>0$, a closed orbit surrounding $P_{0}$ but not $(0,0)$,  \smallskip

\noindent (iii) or an homoclinic orbit defined on $\mathbb{R}$, starting from $(0,0)$ with initial slope
$$\lim_{t\to-\infty}\myfrac{w'(t)}{w(t)}=m$$
where $m$ is defined $E(m)=d$, and ending at
$(0,0)$ with 
$$\lim_{t\to\infty}\myfrac{w'(t)}{w(t)}=-m.$$
\es
\Proof We recall that $E$ and $u$ are defined by (\ref{polu}) and (\ref{D}), by using polar coordinates $(\gr,\gth)$ in the $(w,y)$-plane.

First look at the vector field on the boundary of $\mathcal{Q}$. At any point
$(0,\sigma)$ with $\sigma>0,$ it is given by $(\sigma,0),$ thus it is
transverse and inward. At any point $(\bar{w},0)$ with $\bar{w}>0,$ it is
given by $\left(  0,\bar{w}(b+d-h(\bar{w})\right)  .$ Thus it is transverse
and outward whenever $b+d\leq0$ or $b+d>0$ and $\bar{w}>a,$ and inward
whenever $b+d>0$ and $\bar{w}<a.$

Consider any solution $(w,y)$ of the system, such that $P=(w(0),y(0))\in
\mathcal{Q}$, and let $\left(  \tau_{1},\tau_{2}\right)  $ be its maximal
interval existence in $\mathcal{Q}\mathbf{.}$ At any point $\tau$ where
$u^{\prime}(\tau)=0$ and $u(\gt)>0$, there holds $u^{\prime\prime}(\tau)=-h^{\prime
}(w)w\varphi(u)<0$ from (\ref{suw}). Thus if $\tau$ exists, it is unique, and
it is a maximum for $u$.

Since $w^{\prime}=y>0$, $w$ has the limits $\ell_{2}\in\left(
0,\infty\right] $ as $\gt\uparrow\tau_{2}$ and $\ell_{1}\in\left[  0,\infty\right)$ as $\gt\downarrow\tau_{1}$. Therefore $u$ is strictly monotonous near $\tau_{1},$ and $\tau_{2}$ thus it
has limits $u_{1},u_{2}\in\left[  0,\infty\right],$ in other words $\theta$
has limits $\theta_{1},\theta_{2}\in\left[  0,\pi/2\right]  $

(i) Let us go forward in time. On any interval where $u$ is increasing, one
has $E(\varphi(u))\leq d,$ thus $u$ is bounded and, consequently, $u_{2}$ is finite. If
$\ell_{2}=\infty,$ then $\theta^{\prime}(\gt)\to-\infty,$ as $\gt\uparrow\gt_{2}$; by (\ref{polu}), $\rho$ is
decreasing, thus it is bounded, which is contradictory; thus $\ell_{2}$ is
finite. If $u_{2}>0$ then $\left(  \ell_{2},\ell_{2}\varphi\left(
u_{2}\right)  \right)  $ is stationary, which is impossible. Thus $u$ is
decreasing to $0,$ and the trajectory converges to $(\ell_{2},0).$ If $b+d>0$
and $\ell_{2}=a$, $u^{\prime}$ tends to $0$ from (\ref{suw}), and
\[
u^{\prime\prime}=-(E\circ\varphi)^{\prime}(u)u^{\prime}-h^{\prime}%
(w)w\varphi(u)=-h^{\prime}(a)a\varphi(u)(1+o(1));
\]
therefore $u^{\prime\prime}<0$ near $\tau_{2},$ which is impossible. Finally, either
$b+d\leq0,$ or $b+d>0$ and $\bar{w}>a,$ and $\tau_{2}$ is finite, the
trajectory leaves $\mathcal{Q}$ transversally at $\tau_{2}.$\smallskip 

(ii) Next let us go backward in time. \medskip

$\bullet$ Suppose $u_{1}=0$. Clearly the trajectory converges to $(\ell
_{1},0);$ then necessarily $b+d>0$ and $\ell_{1}\leq a,$ thus $\ell_{1}<a$ as
above. The trajectory enters $\mathcal{Q}$ transversally at $\tau_{2},$ and
from the symmetries it is a closed orbit surrounding only the stationnary point $P_{0}.$

$\bullet$ Next, suppose $u_{1}=\infty$. It means that $\theta$ tends to $\pi/2.$
Then from (\ref{pol}), $\theta^{\prime}$ tends to $1,$ thus $\tau_{1}$ is
finite,
$$\pi/2-\theta=(\tau-\tau_{1})(1+o(1))\;,\;\tan\theta=(\tau-\tau
_{1})^{-1}(1+o(1)),$$ and
\[
(p-1)(\tau-\tau_{1})^{-1}\frac{\rho^{\prime}}{\rho}=(b+1+(d-h(\rho\cos
\theta))\cos^{p-2}\theta)(1+o(1).
\]
If $p\geq2,$ then $\rho^{\prime}/\rho=O((\tau-\tau_{1}));$ if $p<2$ then
$\rho^{\prime}/\rho=O((\tau-\tau_{1})^{p-1}).$ In any case, $\ln\rho$ case is
integrable, thus $\rho$ has a finite limit $\bar{y}>0.$ Then the trajectory
enters $\mathcal{Q}$ transversally at $\tau_{1}$ and from the symmeries it is
a closed orbit surrounding $(0,0).$ From the considerations in \S\ 3-2, for any $\bar{y}>0$
there exist such an orbit, and it is unique. Moreover in $\mathcal{Q}$ the
slope $w'/w=\xi=\varphi(u)$ is decreasing from $\infty$ to $0;$ indeed it decreases
near $\tau_{1}$ and $\tau_{2}$ and can only have a maximal point.

$\bullet$ At end, suppose $0<u_{1}<\infty.$ If $\ell_{1}>0,$ then $\left(
\ell_{1},\ell_{1}\varphi\left(  u_{1}\right)  \right)  $ is stationary, which
is impossible. Thus $(y,w)$ converges to $(0,0).$ And $w^{\prime}/w$ tends to
$\varphi\left(  u_{1}\right),$ thus $\tau_{1}=-\infty.$ And $u^{\prime}$
converges to $d-E(\varphi(u_{1})),$ thus $\tan\theta=\varphi(u)$ has a limit
$m\geq0$ such that $E(m)=d.$ From the symmetries the trajectory is homoclinic
and the solution $w$ is defined on $\mathbb{R}.$\qeda\\

The next theorem studies the precise behaviour of solutions according to the sign of $b+d.$

\bth{Thr1} Assume $p>1$ and consider system (\ref{FG}) in the $(w,y)$-plane.\smallskip

\noindent(i) Assume $b+d>0.$ Then there exists a unique homoclinic trajectory
$\mathcal{H}$ starting from $(0,0)$ in $\mathcal{Q}$ with initial slope
$m_{d}=E^{-1}(d)$ ( $m_{0}=\sqrt{b/(p-1)}$ if $d=0),$ ending at $(0,0)$
with the slope $-m_{d},$ and surrounding $P_{0}.$ Up to the stationary points,
the other orbits are closed, and either they surround only one of the points
$P_{0}$ or $-P_{0},$ in the domain delimitated by $\mathcal{H}$, corresponding
to solutions $w$ of constant sign, or they are exterior to $\pm\mathcal{H}$
and surround $(0,0)$ and $\pm P_{0},$ corresponding to sign changing solutions
$w$.\smallskip

\noindent(ii) Assume $b+d\leq0.$ Then

$\bullet$ if $(p-2)b\leq2(p-1),$ or $\left[  (p-2)b>2(p-1)\text{ and }%
d<E(\eta)\right],$ there is no homoclinic trajectory.

$\bullet$ if $\left[  (p-2)b>2(p-1)\text{ and }E(\eta)<d\leq-b\right]  ;$ then
denoting by $m_{1,d}<m_{2,d}$ the two positive roots of equation $E(m)=d$, there exist infinitely many homoclinic trajectories $\mathcal{H}_{1}$ starting from $(0,0$
in $\mathcal{Q}$ with the initial slope $m_{1,d}$ and ending at $(0,0)$ with the final slope
$-m_{1,d},$ and a unique homoclinic trajectory $\mathcal{H}_{2}$ starting from
$(0,0)$ in $\mathcal{Q}$ with initial slope $m_{2,d}$ and ending at $(0,0)$ with final
slope $-m_{2,d}.$
\es
\Proof
(i) \textbf{Case }$b+d>0.$\textbf{ }Then the equation $E(m)=d$ has a unique
positive solution $m=E^{-1}(d)$; and $w^{\prime}/w$ tends to $m;$ thus the trajectory
starts from $(0,0)$ with a slope $m.$ Then for any $P\in\mathcal{Q},$ the
trajectory $\mathcal{T}_{\left[  P\right]  }$ passing through $P$ meets the
axis $y=0$ after $P$ at some point $(\mu,0)$ with $\mu>a.$ Denote
\begin{equation}\label{UN}
\mathcal{U}\mathbf{=}\left\{  P\in\mathcal{Q}:\mathcal{T}_{\left[  P\right]
}\cap\left\{  (0,\sigma):\sigma>0\right\}  \neq\emptyset\right\}
,\;\mathcal{V}\mathbf{=}\left\{  P\in\mathcal{Q}:\mathcal{T}_{\left[
P\right]  }\cap\left\{  (\mu,0):0<\mu<a\right\}  \neq\emptyset\right\}  .
\end {equation}
Then either $P\in\mathcal{U}$ and the trajectory is a closed orbit surrounding
$(0,0)$ and $\pm P_{0},$ and in $\mathcal{Q}.$ Or $P\in\mathcal{V}$ and the
trajectory is a closed orbit surrounding only $P_{0}$. Or $\mathcal{T}%
_{\left[  P\right]  }$ is an homoclinic orbit $\mathcal{H}$ starting from
(0,0) with the slope $m$, where $m$ is the unique solution of equation
$E(m)=d$ (such that $m>\eta$ if $E$ is not monotone, see (\ref{DerivE})).
Next $\mathcal{U}$ and $\mathcal{V}$ are open, since the vector field  is transverse on the axes, thus $\mathcal{U\cup V}\neq\mathcal{Q}.$ This shows the
existence of such an orbit $\mathcal{H}.\medskip$
 
(ii) \textbf{Case }$b+d\leq0.\medskip$

$\bullet$ Either $b+d<0$ and $E$ is increasing, or $E$ has a minimum at $\eta$
and $d<E(\eta)$. In such a case equation $E(m)=d$ has no solution, and there is no
homoclinic orbit$.$ Or $E$ is increasing and $b+d=0;$ then $E(\varphi
(u))>-b=d,$ thus $u^{\prime}<0,$ thus $u$ cannot tend to 0, and the same
conclusion holds.$\medskip$

$\bullet$ 0r $E$ has a minimum at $\eta$ and $E(\eta)<d\leq-b$. In that case the
equation $E(m)=d$ has two roots $m_{1},m_{2}$ such that $0\leq m_{1}%
<\eta<m_{2}\leq m_{b}$, where $m(b)$ is defined by $E(m_{b})=-b$.  Any trajectory $\mathcal{T}_{\left[  P\right]  }$ such that
$P\in\mathcal{U}$ (see (\ref{UN}) for the definition) satisfies $u^{\prime}<0,$ it means $h(w)>d-E(\varphi(u))$
and the range of $u$ is $(0,\infty),$ therefore there exists $\tau$ such that
$\varphi(u)(\tau)=\eta,$ hence $h(w(\tau))>d-E(\eta)$ and $y(\tau)=\eta
w(\tau).$ Next consider any trajectory $\mathcal{T}_{\left[  \tilde{P}\right]
}$ starting from $\tilde{P}=(\tilde{w},\eta\tilde{w})$ such that $h(\tilde
{w})\leq d-E(\eta)$. It cannot be a trajectory of the preceding
type$,$ thus $(y,w)\rightarrow(0,0)$ as $\tau\rightarrow\tau_{1},$ and
$\theta$ tends to $\theta_{1},$ with $\tan\theta_{1}=m_{1}$ or $m_{2};$
moreover $u^{\prime}(0)\geq0$, and $u^{\prime}<0$ near $\tau_{2},$ thus there
exists a unique $\tau\geq0$ such that $u^{\prime}(\tau)=0;$ then $u^{\prime
}>0$ in $\left(  \tau_{1},\tau\right)  $, therefore $\tan\theta_{1}<\eta,$ and finally
$\tan\theta_{1}=m_{1}.$ Consequently there exist infinitely many such trajectories
$\mathcal{H}_{1},$ with initial slope $m_{1}$. Next fix one trajectory
$\mathcal{T}_{\left[  \tilde{P}_{0}\right]  }$ such that $h(\tilde{w}_{0})\leq
d-E(\eta)$. Let $\mathcal{R}$ be the subdomain of $\mathcal{Q}$ delimitated by
$\mathcal{T}_{\left[  \tilde{P}_{0}\right]  }$ and $\mathcal{T}_{\left[
(0,1)\right]  }$ and
\[
\mathcal{V}\mathbf{=}\left\{  P\in\mathcal{R}:\mathcal{T}_{\left[  P\right]
}\cap\left\{  (w,\eta w):0<w<\tilde{w}_{0}\right\}  \neq\emptyset\right\}  .
\]
The set $\mathcal{V}$ is open because the intersection with the line $y=\eta
w$ for $w<$ $\tilde{w}_{0}$ is transverse since at the intersection point,
$h(\tilde{w})<d-E(\eta),$ thus $u^{\prime}>0,$ and $y/w=\varphi(u)=\eta,$ and
\[
\frac{y^{\prime}}{y}=\varphi(u)+\frac{\varphi^{\prime}(u)}{\varphi
(u)}u^{\prime}>\eta=\frac{w^{\prime}}{w}.%
\]
Then ($\mathcal{U\cap R)\cup V}\neq\mathcal{R}.$ Then there exists at least a
trajectory $\mathcal{H}_{1,\ast}$ starting from $(0,0)$ with initial slope
$m_{2}.\medskip$

(iii) \textbf{Uniqueness of} $\mathcal{H}$ and $\mathcal{H}_{2}$. Let
$m=m_{0}$ or $m_{2,d}.$ Suppose that system (\ref{FG}) has two solutions
$(w_{1},y_{1}),$ $(w_{2},y_{2})$ defined near $-\infty,$ such that $w_{i}>0$ and
$w_{i}(\gt)$ tends to $0$ and $y_{i}(\gt)/w_{i}(\gt)$ tend to $m$ as $\gt\downarrow -\infty$. Then the system
(\ref{suw}) has two local solutions $(w_{1},u_{1}),$ $(w_{2},u_{2})$ such that
$\varphi(u_{i})$ tends to $m$ at $-\infty.$ Then $w_{i}^{\prime}>0$ locally and one can express $u_{i}$ as a function of $w_{i}.$ Then at the same point
$w,$
\[
w\frac{d(\Psi(u_{i})}{dw}=w\varphi(u_{i})\frac{du_{i}}{dw}=-E(\varphi
(u_{i}))-h(w)+d,
\]%
\[
w\frac{d(\Psi(u_{2})-\Psi(u_{1}))}{dw}=-(E(\varphi(u_{2}))-E(\varphi
(u_{1}))=E^{\prime}(\varphi(u^{\ast}))\varphi^{\prime}(u^{\ast})(u_{2}-u_{1})
\]
for some $u^{\ast}$ between $u_{1}$ and $u_{2},$ and $E^{\prime}%
(\varphi(u^{\ast}))=E^{\prime}(m)(1+o(1));$ and $E^{\prime}(m)>0.$ Then for
small $w$
\[
\frac{d(\Psi(u_{2})-\Psi(u_{1}))}{dw}(\Psi(u_{2})-\Psi(u_{1}))<0,
\]
which implies that ($\Psi(u_{2})-\Psi(u_{1}))^{2}$ is decreasing, with limit $0$ at $0$.
Therefore $\Psi(u_{2})=\Psi(u_{1}),$ thus $u_{2}\equiv u_{1}$ near $-\infty$; but
from (\ref{suw}), $h(w_{1})=h(w_{2}),$ and since $h$ is one to one, it follows$w_{1}\equiv
w_{2}$ near $-\infty$. The global uniqueness follows, since the system is
regular except at $(0,0).$ All the trajectories are described.\qeda\\

\noindent\Remark
Under the assumption (\ref{lfp}), existence and uniqueness of $\mathcal{H}$
and $\mathcal{H}_{2}$ can be obtained in a more direct way whenever $d\neq E(\eta).$
Indeed the system (\ref{suvt}) relative to $(v,u)$ is regular, with stationary
points $(0,0),$ $(0,\pm\varphi^{-1}(m)),$ where $m=m_{0},m_{1}$ or $m_{2}$ and
also $(\pm a,0)$ if $b+d>0.$ The linearized system at $(0,\varphi^{-1}(m))$ is
given by the matrix $%
\begin{pmatrix}
m(q+1-p) & 0\\
0 & K(m)
\end{pmatrix}
,$ with $K(m)=p(p-1)(\eta^{2}-m^{2})/(1+(p-1)m^{2}).$ If $m=m_{1,d},$ then it
is a source, and we find again the existence of an infinity of solutions. If
$m=m_{d}$ or $m=m_{2,d},$ then $K(m)<0,$ thus this point is a saddle point.
Then in the phase plane $(v,u),$ there exists precisely one trajectory defined
near $-\infty,$ such that $v>0$ and converging to $(0,m)$ at $-\infty,$ and
$u/v$ converges to $0.$\qeda\\

\noindent\Remark
Suppose $f(w)=|w|^{q-1}w,$ then we can study the critical case $(p-2)b>2(p-1)$
and $E(\eta)=d:$ there exist infinitely many homoclinic trajectories
$\mathcal{H}_{1}$ starting from $(0,0)$ in $\mathcal{Q}$ with an infinite initial
slope and ending at $(0,0)$ with an infinite slope, and a unique homoclinic
trajectory $\mathcal{H}_{2}$ starting from $(0,0)$ in $\mathcal{Q}$ with the initial
slope $\eta$ and ending at $(0,0)$ with the slope $-\eta.$ Indeed using system
(\ref{suv}) and setting $u=\varphi^{-1}(\eta)+z,$ and $\zeta=(q+1-p)\eta z+v,$
it can be written under the form%
\[
\zeta^{\prime}=P(\zeta,v),\qquad\qquad v^{\prime}=(q+1-p)\eta v+Q(\zeta,v),
\]
where $P$ and $Q$ both start with quadratic terms. Moreover the quadratic part
of $P(\zeta,v)$ is given by $p_{2,0}\zeta^{2}+p_{1,1}\zeta v+p_{0,2}v^{2},$
where by computation,
\[
p_{2,0}=-\frac{p(p-1)}{q+1-p}\eta\varphi^{\prime2}(\varphi(\eta))(1+\eta
^{2})^{(p-2)/2}<0.
\]
The results follow from the description of sadle-node behaviour given in
\cite[Theorem 9.1.7]{HuWe}.$\medskip$

\noindent\Remark
In the case $b=1>-d,$ we have a representation of the homoclinic trajectory :
it corresponds to $K=0$ in (\ref{emo})$.$ In the case $f(w)=|w|^{q-1}w$, in
terms of $u$ we obtain%
\[
u^{\prime}=\frac{q+1-p}{p}(E(\varphi(u))-d),
\]
which allows to compute $u$ by a quadrature.

\subsection{Period of the solutions}

First we consider the sign changing solutions

\bth{ts} Assume $p>1$. For any $\nu>0$ let $\mathcal{T}_{\left[  (0,\nu)\right]  }$  be the trajectory which starts from $(0,\nu)$, and let  $T(\nu)$ be its least period. Then $\gn\mapsto T(\gn)$ is decreasing on $\left(  0,\infty\right).$ Furthermore the range of $T(.)$ can be computed in the following way.\medskip

\noindent(i) If $b+d\leq 0$ and $m\mapsto E (m)$ is increasing, or if $d<\min E,$ then $T(.)$ decreases
from $T_{d}$ to $0$, where%
\begin{equation}
T_{d}=4%
{\displaystyle\int\nolimits_{0}^{\infty}}
\frac{du}{E(\varphi(u))-d}=4%
{\displaystyle\int\nolimits_{0}^{\pi/2}}
\frac{1+(p-1)\tan^{2}\theta}{(p-1)\tan^{2}\theta-b-d\cos^{p-2}\theta}d\theta,
\label{td}%
\end{equation}
and $T_{d}$ is finite if and only if $b+d<0$. If $b<0=d,$ then
\begin{equation}
T_{0}=2\pi\frac{(p-1)\gamma+1}{(p-1)\gamma(\gamma+1)}\;\text { with }\;\gamma
=\sqrt{\left\vert b\right\vert /(p-1)}). \label{to}%
\end{equation}
(ii) If $b+d>0$ or $b+d\leq0$ and $d\geq\min E,$ then $T(.)$ decreases from $\infty$ to
$0$.
\es
\Proof
{\it{Step 1. Monotonicity of }$T$}. Consider the part of the trajectories
$\mathcal{T}_{\left[  (0,\nu)\right]  }$ located in $\mathcal{Q},$ given by
$(w_{\nu},y_{\nu}).$ We have already shown that $u$ is decreasing with respect to
$\tau$ from $\infty$ to $0,$ then $E(\varphi(u))+h(w_{\nu}(u))-d>0$ and
$w_{\nu}$ can be expressed in terms of $u,$ and
\begin{equation}\label{T-E}
T(\nu)=4%
{\displaystyle\int\nolimits_{0}^{\infty}}
\frac{du}{E(\varphi(u))+h(w_{\nu}(u))-d}.%
\end{equation}
Let $\lambda>1.$ Since the trajectories $\mathcal{T}_{\left[  (0,\nu)\right]
}\,$and $\mathcal{T}_{\left[  (0,\lambda\nu)\right]}$ have no intersection
point, $w_{\lambda\nu}(u)>w_{\nu}(u)$ for any $u>0,$ and $h$ is
nondecreasing, thus $T(\lambda\nu)<T(\nu),$ and $T$ is decreasing.\medskip

\noindent{\it{Step 2. Behaviour near }$\infty$}.  Let $\nu_{n}\geq1,$ such that
$\lim\nu_{n}=\infty.$ Observe that for fixed $u,$ for any integer $n\geq1,$
there exists a unique $\tilde{\nu}_{n}>0$ (depending on $u),$ such that
$w_{\tilde{\nu}_{n}}(u)=n;$ let $\hat{\nu}_{n}=\max(\tilde{\nu}_{n},n).$ Then
$h(w_{\hat{\nu}_{n}}(u))\geq h(n)$, thus $h(w_{\hat{\nu}_{n}}(u))$ converges to
$\infty;$ since $\nu\mapsto h(w_{\nu}(u))$ is nondecreasing then $h(w_{\nu
_{n}}(u))$ converges to $\infty,$ and $T(\nu_{n})$ converges to $0$, using the
Beppo-Levi theorem.\medskip

\noindent{\it Step 3. Behaviour near $0.$}

$\bullet$ First assume $b+d\leq0,$ and $E$ is increasing, or $d<E(\eta).$ Then
all the orbits are of the type $\mathcal{T}_{\left[  (0,\nu)\right]}.$ Let
$\nu_{n}\in\left(  0,1\right),$ such that $\lim\nu_{n}=0.$ For fixed $u$
and any integer $n\geq1,$ there exists a unique $\bar{\nu}_{n}>0$ (depending
on $u),$ such that $w_{\bar{\nu}_{n}}(u)=1/n;$ let $\check{\nu}_{n}=\min
(\bar{\nu}_{n},1/n).$ Then $h(w_{\check{\nu}_{n}}(u)\leq h(1/n)$, thus
$h(w_{\check{\nu}_{n}}(u))$ converges to $0,$ and again $h(w_{\nu_{n}}(u))$
converges to $0.$ Then $T(\nu_{n})$ converges to $T_{d}$ given by (\ref{td}),
using the Beppo-Levi theorem. If $b+d<0,$ then $T_{d}$ is finite: indeed near
$\infty,$ $E(\varphi(u))=(p-1)u^{p/(p-1)}(1+o(1));$ if $E$ is increasing, then
$E(\varphi(u))-d>-(b+d)>0;$ if $d<E(\eta),$ then $E(\varphi(u))-d\geq
E(\eta)-d>0$.

If $b+d=0$ and $E$ is increasing, then $T_{d}=\infty$: indeed near $0,$
$$E(\varphi(u))-d=u^{2}(E^{\prime\prime}(0)/2+o(1))$$ 
and 
$$E^{\prime\prime}(0)=2(p-1)-(p-2)b\;\text{if }\;(p-2)b=2(p-1).$$ 
Therefore 
$$E(\varphi(u))-d=(p(p-1)/4)u^{4}(1+o(1)).$$
In all the cases the integral (\ref{T-E}) giving $T$ is divergent.

When $b<0=d,$ one can compute $T_{0}:$%
\begin{align*}
\frac{T_{0}}{4}  &  =%
{\displaystyle\int\nolimits_{0}^{\infty}}
\frac{du}{E(\varphi(u))}=%
{\displaystyle\int\nolimits_{0}^{\infty}}
\frac{\phi^{\prime}\left(  \xi\right)  d\xi}{E(\xi)}=%
{\displaystyle\int\nolimits_{0}^{\infty}}
\frac{1+(p-1)\xi^{2}}{(\left\vert b\right\vert +(p-1)\xi^{2})(1+\xi^{2})}%
d\xi\\
&  =\frac{\pi}{2}+(\frac{1}{p-1}-\gamma^{2})%
{\displaystyle\int\nolimits_{0}^{\infty}}
\frac{ds}{(\gamma^{2}+s^{2})(1+s^{2})}=\frac{\pi}{2}(1+\frac{1-(p-1)\gamma
^{2}}{(p-1)\gamma(\gamma+1)}).
\end{align*}
Hence (\ref{to}) holds.

$\bullet$ Next assume $d>E(\eta).$ Considering $\nu_{n}$ as above, for any
fixed $u$ such that $\varphi(u)>m_{2},$ there exists a unique $\bar{\nu}%
_{n}>0$ (depending on $u),$ such that $w_{\bar{\nu}_{n}}(u)=1/n.$ As
above,
\[%
{\displaystyle\int\nolimits_{\varphi^{-1}(m_{2})}^{\infty}}
\frac{du}{E(\varphi(u))+h(w_{\bar{\nu}_{n}})-d}\rightarrow%
{\displaystyle\int\nolimits_{\varphi^{-1}(m_{2})}^{\infty}}
\frac{du}{E(\varphi(u))-d}=\infty,%
\]
since  $E'(m_{2})$ is finite. As a consequence, $T(\nu_{n})$ tends to $\infty.$ If $d=E(\eta),$ the
same proof still works with $m_{2}$ replaced by $\eta:$ the integral is still
divergent because the denominator is of order $2$ in $u-\varphi^{-1}(\eta)$, as,  near
$0,$ there holds
$$E(\varphi(u))-d=\frac{1}{2}E^{\prime\prime}\left(  \eta\right)
(\varphi(u)-\eta)^{2}(1+o(1))=\frac{1}{2}E^{\prime\prime}\left(  \eta\right)
(\varphi(u)-\eta)^{2}(1+o(1),$$
and $$E^{\prime\prime}\left(  \eta\right)
=2p(p-1)\eta^{2}(1+\eta^{2})^{(p-2)/2}>0.$$ 
At last suppose $b+d>0;$ the same
proof with $m_{2}$ replaced by $m$ shows that $T(\gn)$ converges to
$\infty$ as $\nu$ tends to $0,$ since $E'(m)$ at $m=E^{-1}(d)$ is
finite.\qeda\bigskip

The monotonicity of the period function is a more general
property, since we have the following result.

\bprop{Lem1}
Let $F,G\in C^{1}(\mathbb{R}^{2}\backslash(0,0))$  are such that $F$ (resp.
$G$) is odd with respect to $y$ (resp. $x$) and even with respect to $x$
(resp. $y),$ with $F(w,y)>0$ in $\mathcal{Q}.$ Assume that for any
$(w,y)\in\mathcal{Q},$ and any $\lambda>0,$%
\begin{equation}
\frac{\partial}{\partial\lambda}\left(  \frac{F(\lambda w,\lambda y)}{\lambda
}\right)  \geq0\text{ }(\text{resp.}\leq 0)\quad\text{and }\frac{\partial
}{\partial\lambda}\left(  \frac{G(\lambda w,\lambda y)}{\lambda}\right)
<0\; (\text{resp.}>0). \label{hyp}%
\end{equation}
Assume also that for any $\sigma$ in some interval $\left(
\sigma_{1},\sigma_{2}\right)$ (where $0<\sigma_{1}<\sigma_{2}$), the trajectory
$\mathcal{T}_{\left[  (0,\sigma)\right]  }$ of solution of the system
\begin{equation}\label{exten}\left\{\BA {l}
w^{\prime}=F(w,y)\\
y^{\prime}=G(w,y)
\EA\right.\end{equation}
passing through $(0,\sigma)$
(necessarily entering $\mathcal{Q}$ since $F(0,\sigma)>0)$ leaves
$\mathcal{Q}$ transversally in a finite time $T(\sigma)/4$ at some point
$(c(\sigma),0)$ (thus $G(c(\sigma),0)<0).$ Then (from the symmetries),
$\mathcal{T}_{\left[  (0,\sigma)\right]  }$ is a closed orbit surrounding
$(0,0),$ with period $T(\sigma),$ and $\gs\mapsto T(\gs)$ is decreasing (resp. increasing) on $\left(  \sigma_{1},\sigma_{2}\right).$
\es

\noindent\Remark We can notice the condition on $F$ is equivalent to $F(\lambda w,\lambda y)\geq\lambda
F(w,y)$ for any $\lambda>1$. The second condition implies that for any
$\lambda>1,$
\[
G(\lambda w,\lambda y)<\lambda G(w,y)\text{ (resp.}G(\lambda w,\lambda
y)>\lambda G(w,y)).
\]
\medskip

\noindent{\it Proof of \rprop {Lem1}}.
In polar coordinates $(\rho,\theta)$ in $\mathcal{Q},$ we get
\[
\rho^{\prime}=F\cos\theta+G\sin\theta,\qquad\qquad\theta^{\prime}=\frac
{1}{\rho}(G\cos\theta-F\sin\theta).
\]
At each point $\tau$ where $\theta^{\prime}(\tau)=0,$ there holds%
\[
\rho\theta^{\prime\prime}(\tau)=\left(  \frac{\partial G}{\partial\rho}%
\cos\theta-\frac{\partial F}{\partial\rho}\sin\theta\right)\;,\;\;  \rho^{\prime
}(\tau)=\frac{F}{\cos\theta}\left(  \frac{\partial G}{\partial\rho}\cos
\theta-\frac{\partial F}{\partial\rho}\sin\theta\right).
\]
But (\ref{hyp}) is equivalent to $\partial F/\partial\rho\geq F/\rho$ and
$\partial G/\partial\rho<G/\rho$ (resp $>$), thus%
\[
\rho\theta^{\prime\prime}(\tau)<\frac{F}{\rho\cos\theta}(G\cos\theta
-F\sin\theta)=0\quad(\text{resp.}>)
\]
In both case $\theta^{\prime\prime}$ has a constant sign. But $\theta^{\prime
}(0)=-F(0,\sigma)<0$ and $\theta^{\prime}(\sigma)=G(c(\sigma),0)<0$ thus
 we get a contradiction by considering the first (resp. the last)
point where $\theta^{\prime}(\tau)=0,$ which satisfies $\theta^{\prime\prime
}(\tau)\geq0$ (resp. $\leq0).$ Thus $\theta$ is decreasing from $\pi/2$ to
$0.$ Then the curves can be represented in function of $\theta$ by $\left(
\rho\left(  \sigma,\theta\right) ,\theta(\sigma)\right),$ and
\[
T(\sigma)=4%
{\displaystyle\int\nolimits_{0}^{\pi/2}}
\frac{d\theta}{H(\rho\left(  \sigma,\theta\right) ,\theta)}%
\]
with
\[
H(\rho,\theta)=\frac{1}{\rho}(F(\rho\cos\theta,\rho\sin\theta)\sin
\theta-G(\rho\cos\theta,\rho\sin\theta)\cos\theta)
\]
Let $\lambda>1.$ Since the trajectories $\mathcal{T}_{\left[  (0,\sigma
)\right]  }\,$and $\mathcal{T}_{\left[  (0,\lambda\sigma)\right]  .}$ have no
intersection point, then $\rho\left(  \lambda\sigma,\theta\right)
>\rho\left(  \sigma,\theta\right)  $ for any $\theta\in\left(  0,\pi/2\right)
;$ by assumption, for fixed $\theta,$ the function $\rho\mapsto F(\rho
\cos\theta,\rho\sin\theta)/\rho$ is nondecreasing (resp. nonincreasing) and
$\rho\mapsto G(\rho\cos\theta,\rho\sin\theta)/\rho$ is decreasing (resp.
increasing), thus $H(\rho\left(  \lambda\sigma,\theta\right) ,\theta
)>H(\rho\left(  \sigma,\theta\right) ,\theta)$, which yields to $T(\lambda
\sigma)<T(\sigma)$ (resp. $>)$. This implies that $T$ is decreasing (resp. increasing).\qeda\\

Next we consider the positive solutions $w$ of equation (\ref{E}).

\bprop{lit} Assume $p>1$ and $b+d>0.$ Consider the trajectories $\mathcal{T}_{\left[
(\mu,0)\right]  }$in the phase plane $(w,y)$ which goes through $(\mu,0),$
for some $\mu\in\left(  0,a\right).$ Let $T^{+}(\mu)$ be their least period. Then
\[
\lim_{\mu\rightarrow0}T^{+}(\mu)=\infty,\qquad\lim_{\mu\rightarrow a}T^{+}%
(\mu)=\frac{2\pi}{\sqrt{ah^{\prime}(a)}}.
\]
In particular if $f(w)=|w|^{q-1}w$, then $\lim_{\mu\rightarrow a}T^{+}(\mu)=$
$2\pi/(q+1-p)(b+d)$.
\es
\Proof
We notice that the trajectory $\mathcal{T}_{\left[  (\mu,0)\right]  }$ intersects the line $y=0$ at
$(\mu,0)$ and another point $(g(\mu),0),$ with $\mu<a<g(\mu),$ and $g$ is decreasing.\smallskip

\noindent{\it Step1. Behaviour near }$a$. When $\mu$ tends to $a,$
then also $g(\mu)$ tends to $a.$ Indeed for any small $\varepsilon>0,$ then
$g(\mu)-a<\varepsilon$ as soon as $\mu-a<\min(\varepsilon,a-g^{-1}%
(a+\varepsilon)).$ Since, along such a trajectory in
$\mathcal{Q}$, $\xi=\varphi(u)$ varies from $0$ to 0, it has a maximal  $\xi^{\ast}$, where $u^{\prime}=0,$ thus
$E(\xi^{\ast})=h(w^{\ast}).$ When $\mu$ tends to $a,$ then $h(w^{\ast})$ tends
to $b,$ thus $\xi^{\ast}$ tends to $E^{-1}(b)=0,$ thus also $\max
_{y\in\mathcal{T}_{\left[  (\mu,0)\right]  }}\left\vert y\right\vert $ tends
to $0$. Using the linearized form of the system at $P_{0}$, and  polar
coordinates with center $(a,0),$ $w=a+r\cos\eta,$ $y=\sqrt{ah^{\prime}(a)}r\sin\eta,$
then $r$ tends to $0$ as $\mu$ tends to $a,$ and one finds $\eta^{\prime
}=-\sqrt{ah^{\prime}(a)}+R/r,$ where $R$ involves the derivatives of $G$ of
order 2, which are bounded near the point $(a,0),$ thus $R/r^{2}$ is bounded.
Therefore $\eta^{\prime}$ tends to $-\sqrt{ah^{\prime}(a)},$ and finally $T^{+}(\mu)$
tends to $2\pi/\sqrt{ah^{\prime}(a)}.\medskip$

\noindent{\it Step2. Behaviour near }$0$. On the trajectory $\mathcal{T}_{\left[
(\mu,0)\right]  },$ the function $u$ is increasing up to a maximal value
$u^{\ast}(\mu),$ and then decreasing; moreover $u^{\ast}$ is a nonincreasing
function of $\mu,$ because two different trajectories have no intersection. Let
$\mu_{n}\in\left(  0,a\right),$ such that $\lim\mu_{n}=0.$ For any $n$ there
exists $\tilde{\mu}_{n}\in\left(  0,a\right)  $ such that the orbit $\mathcal{T}_{\left[(\tilde{\mu}_{n},0)\right]  },$ contains a
point above the line $y=\varphi^{-1}(m)(1-1/n)w,$ let $\hat{\mu}_{n}=\min
(\mu_{n},1/n).$ Then $u^{\ast}(\hat{\mu}_{n})\geq\varphi^{-1}(m)(1-1/n),$ thus
$u^{\ast}(\mu_{n})$ tends to $m$ ; then from the Beppo-Levi theorem
\[
\lim\inf T^{+}(\mu)\geq\lim%
{\displaystyle\int\nolimits_{u^{\ast}(\mu)}^{\infty}}
\frac{du}{E(\varphi(u))-d+h(w(\mu,u))}=%
{\displaystyle\int\nolimits_{m}^{\infty}}
\frac{du}{E(\varphi(u))-d+h(w(u))}%
\]
where $w$ is the solution defining $\mathcal{H},$ and this integral is infinite.
\qeda\\

\noindent\Remark
Here the question of the monotonicity of the period is difficult to answer, even for
$p=2$, where it is solved by using the first integral, see \cite{BiBo}. It is
open in the general case. More generally, if a dynamical system a center, the description of the period function is still a chalenging problem. For example, one can contruct a quadratic dynamical system with a center, the
associated period function of which is not monotone, and even with at least two
critical points, see \cite{ChDu} and \cite{ChJa}.\\

\noindent\Remark
In the case $b=1,$ we can compute theoretically the period $T^{+}$ by using
the first integral (\ref{pain}). The stationary point $P_{0}=(h^{-1}(1),0)$ is
obtained for $K_{a}=a^{p}/p-\mathcal{F}(a)>0$ (in case of a power,
$K_{a}=(q+1-p)/p(q+1)).$ The positive solutions correspond to trajectories
$\mathcal{T}_{K}$ with $K\in\left(  0,K_{a}\right),$ intersecting the axis $y=0$
at points $\left(  w_{1},0\right)  $, $\left(  w_{2},0\right)  $ with
$w_{1}<a<w_{2}$ defined by $w_{i}^{p}/p-\mathcal{F}(w_{i})=K,$ and the period
is given by
\[
T^{+}=2%
{\displaystyle\int\nolimits_{w_{1}}^{w_{2}}}
\frac{dw}{wE^{-1}(-p\frac{K+\mathcal{F}(w)}{w^{p}})}.
\]
Unfortunately, this formula does not allow us to prove the monotonicity of the period function for $p\neq2.$\\

It is remarkable that, in the case $f(w)=\left\vert w\right\vert ^{q-1}w,$ one can solve completely the problem
in the particular case where $b=1$ and $q=2p-1,$ using the equation (\ref{Eu})
satisfied by $u$.

\bprop{Prop2}
Suppose that $f(w)=\left\vert w\right\vert ^{q-1}w,$ and $b=1$ and $q=2p-1,$
$p>1$ and $d+1>0.$ If $p>2$ or $d+1<1/(2-p),$ then $T^{+}$ is decreasing on
$(0,a).$
\es
\Proof
Since $B(\xi)=0$ by (\ref{bks}), equation (\ref{Eu}) turns to
\[
u^{\prime\prime}=(q+1-p)(E(\varphi(u)-d)\varphi(u)=(q+1-p)(-(1+d)+p\Psi.
(u))\varphi(u)
\]
Henceforth
\[
\frac{1}{q+1-p}u^{\prime\prime}u^{\prime}=-(1+d)\Psi^{\prime}(u)u^{\prime
}+p\Psi(u)\Psi^{\prime}(u)u^{\prime},%
\]
from which expression we derive the first integral,
\begin{equation}
\frac{1}{q+1-p}u^{\prime2}=C-\mathcal{U}(u)),\qquad\mathcal{U=\mathcal{M}%
}\circ\Psi,\qquad\mathcal{M}(t)=2(1+d)t-pt^{2}. \label{upr}%
\end{equation}
From (\ref{upr}) the integral curves $\mathcal{S}$ in the $(u,u')$-plane are symmetric with respect to the
axis $u^{\prime}=0.$ The times for going from $u=0$ to $u=u^{\ast}$ and from
$u^{\ast}$ to $0$ are equal, and $u^{\ast}$ is given by $C=\mathcal{M}%
(\Psi(u^{\ast})).$ The computation of the period is reduced to the part relative
to the first quadrant. Here we follow the method of \cite{BiBo}: we get
\[
T^{+}(u^{\ast})=4%
{\displaystyle\int\nolimits_{0}^{u^{\ast}}}
\frac{d\eta}{\sqrt{\mathcal{U}(u^{\ast})-\mathcal{U}(\eta)}}=4%
{\displaystyle\int\nolimits_{0}^{1}}
\frac{u^{\ast}ds}{\sqrt{\mathcal{U}(u^{\ast})-\mathcal{U}(su^{\ast})}}.
\]
Then
\[
\frac{dT^{+}(u^{\ast})}{du^{\ast}}=4%
{\displaystyle\int\nolimits_{0}^{1}}
\frac{(\Theta(u^{\ast})-\Theta(su^{\ast}))ds}{(\mathcal{U}(u^{\ast
})-\mathcal{U}(su^{\ast}))^{3/2}},\qquad\text{with\quad\ }\Theta(u^{\ast
})=\mathcal{U}(u^{\ast})-u^{\ast}\mathcal{U}^{\prime}(u^{\ast})/2,
\]
and
\[
2\frac{d\Theta(u^{\ast})}{du^{\ast}}=2\Theta^{\prime}(u^{\ast})=\mathcal{U}%
^{\prime}(u^{\ast})-u^{\ast}\mathcal{U}^{\prime\prime}(u^{\ast}).
\]
In the interval of study, $\varphi(u^{\ast})<E^{-1}(d)$, $(E\circ
\varphi)(u)<d$ from (\ref{D}), thus $\Psi(u)<(1+d)/p,$ and $\mathcal{M}$ is
increasing for $0<t<(1+d)/p,$ thus $\mathcal{U}^{\prime}>0.$ Then at any point
$u$, $\Theta^{\prime}(u)>0\Longleftrightarrow(\mathcal{U}^{\prime}/u)^{\prime
}<0.$ Now
\[
\frac{\mathcal{U}^{\prime}(u)}{2pu}=\frac{(-E(\varphi(u))+d)\varphi(u)}%
{u}=1-(p-1)\varphi^{2}(u)+d(1+\varphi^{2}(u))^{(2-p)/2},%
\]
hence $(\mathcal{U}^{\prime}/u)^{\prime}=2X(u)\varphi(u)\varphi^{\prime}(u),$
with
\[
X(u)=-(p-1)+(2-p)d(1+\varphi^{2}(u))^{-p/2},%
\]
and $d>E(\varphi(u));$ it implies $X(u)<0$ if $p>2$ or $p<2$ and
$d<(p-1)/(2-p).$ Henceforth $\Theta$ is increasing, and the same holds for $P$ as a function of  $u^{\ast}$. Finally $u^{\ast}$ is decreasing with respect to $\mu,$
and consequently $P$ is decreasing with respect to $\mu.$
\qeda\\

\noindent\Remark
When $p=2,$ and $q=2p-1=3,$ equation (\ref{Eu}) reduces to $u^{\prime\prime
}=-2u+2u^{3},$ which, surprisingly, is an equation correponding to the problem with absorption, and (\ref{E}) reduces to $w^{\prime\prime}-w+w^{3}=0$. In this case, all the solutions can be expressed in terms of elliptic integrals, see \cite{BiBo}.

\subsection{Returning to the initial problem}

\noindent {\it Proof of Theorem 2. } Here $\beta=\beta_{q}=p/(q+1-p)$, 
$\lambda=\lambda_{q}$ is given by (\ref{2dim3}) and $c_{q}=\beta_{q}^{p-2}%
\lambda_{q}$ by (\ref{2dim4}). Moreover $\omega\left(  \sigma\right)  =\beta
_{q}^{\beta_{q}}w(\beta_{q}\sigma)$ from (\ref{cv}), $b=-\lambda_{q}%
/\beta_{q}^{2}=-c_{q}/\beta_{q}^{p}$ and $d=c/\beta_{q}^{p}$ from (\ref{bd}). At end
$f(s)=g(s)=\left\vert s\right\vert ^{q-1}s$ and $h(s)=$ $\left\vert
s\right\vert ^{q-p}s.$ Thus $c>c_{q}$ is equivalent to $b+d>0,$ and then the
constant solutions $w\equiv\pm(b+d)^{1/(q-p+1)}$ of (\ref{E}) correspond to
the constant solutions $\omega\equiv\pm\left(  c-c_{q}\right)  ^{1/(q+1-p)}$
of equation (\ref{2dim2}). For any $\operatorname{integer}$ $k\geq1,$ we look for
periodic solutions $\omega$ of smallest period $2\pi/k,$ or equivalently
solutions $w$ of period $T_{k}=2\pi\beta_{q}/k.$ From (\ref{DerivE}), the
function $E$ is increasing. First consider the sign changing solutions: if
$c\geq c_{q}$, then from \rth{ts}, the period function $T$ of $w$
is decreasing from $\infty$ to 0, hence for any $k\geq1$ it takes precisely
once the value $T_{k}.$ If $c<c_{q},$ then $T$ decreases from $T_{d}$ given by
(\ref{td}) to $0,$ thus it takes once the value $T_{k}$ for any $k>M_{q}%
=T_{d}/2\pi\beta_{q}$ given at (\ref{2dim6}). Next consider the positive
solutions: from \rprop{lit}, the period function of $w$ takes any
value between $\infty$ and $2\pi/\sqrt{(q+1-p)(b+d)},$ thus it takes the value
$T_{k}$ for any $k$ $<(p\beta_{q}^{1-p}(c-c_{q}))^{1/2},$ which ends the proof.
\qeda\\

In the case of equation (\ref{2dim2}) (i.e. $c=0)$, we
obtain the following description of the sets $\CE$ and $\CE^+$:

\bcor{13} Assume $p>1,$ $q>p-1$, and $c=0.$\smallskip

\noindent(i) Then the set $\mathcal{E}$ of changing sign solutions of  (\ref{2dim2}) is given by (\ref{2dim5}), where $k_{q}=1$ if $p<2$ and
$q\geq2(p-1)/(2-p),$ and $k_{q}>M_{q}$ if $p\geq2$ or ($p<2$ and
$q<2(p-1)/(2-p)),$ where 
\begin{equation}M_{q}=2/(q-1),\label{Mq0}%
\end{equation}
if $p=2$, and
\begin{equation}
M_{q}=\frac{(p-2)m_{q}}{((p-1)m_{q}+1)(m_{q}-1)},\;\text{ with }\; m_{q}=\sqrt
{\frac{(2(p-1)+(p-2)q}{p(p-1)}}, \label{Mq}%
\end{equation}
if $p\neq2$.

(ii) If $p\geq2$ or ($p<2$ and $q<2(p-1)/(2-p)),$ then $\mathcal{E}%
^{+}=\varnothing.$ If $p<2$ and $q\geq2(p-1)/(2-p),$ then $\mathcal{E}%
^{+}=\left\{  (-c_{q})^{1/(q+1-p)}\right\}  .$\bigskip
\es
\Proof
Here $c_{q}<0$ is equivalent to $p<2$ and
$q>2(p-1)/(2-p).$ Furthermore $M_{q}=T_{0}/2\pi\beta_{q}$ can be computed from
(\ref{to}), which gives (\ref{Mq0}), (\ref{Mq}). Moreover in any case $c_{q}+$ $\beta
_{q}^{p-1}/p=$ $\beta_{q}^{p}(p-1)(q+1)/p^{2}>0$ thus there exist no positive
nonconstant periodic solutions.\qeda\\

\noindent {\it Proof of Corollary 1. } Let $S$ be a sector on $S^1$ with opening angle $\gth\in (0,2\gp)$. From \cite [Th 3.3]{KV}, $\beta_{S}$ is the positive solution of equation
$$\phi(\beta_{S})=\left(1+\myfrac{1}{k}\right)^2\left(\beta_{S}^2+\myfrac{p-2}{p-1}\beta_{S}-(\beta_{S}+1)^2\right)=0,
$$
where $k=\gp/\gth\geq 1$. Using \rcor{13} (applied without assuming that $k$ is an integer) we distinguish two cases:\smallskip

(i) $p<2$ and $q\geq 2(p-1)/(2-p)$. Then there always exists a solution to the Dirichlet problem in $S$. Notice that $0<\beta_{q}\leq (2-p)/(p-1)$, thus $\phi(\beta_{S})<0$ and consequently $\beta_{q}<\beta_{S}$.\smallskip

(ii) $p>2$ or $p<2$ and $q<2(p-1)/(2-p)$. The existence is equivalent to $k>M_{q}$ 
(see (\ref{Mq})). It means
$$\left(1+\myfrac{1}{k}\right)^2<\left(\myfrac{(p-1)m_{q}^2-1}{m_{q}(p-2)}\right)^2
=\myfrac{(\beta_{q}+1)^2}{\beta_{q}^2m_{q}^2}
=\myfrac{(\beta_{q}+1)^2}{\beta_{q}(\beta_{q}+(p-2)\beta_{q}/(p-1))}.
$$
Thus $\phi(\beta_{q})<0$. Equivalently, $\beta_{q}<\beta_{S}$.\qeda 
\section{The case $p=1$\label{ega}}

\subsection{Existence of a first integral}

As shown in \rlemma{lem}, we can reduce the study to
\begin{equation}
\frac{d}{d\tau}\left(  \frac{w^{\prime}}{\sqrt{w^{2}+w^{\prime2}}}\right)
-\,b\frac{w}{\sqrt{w^{2}+w^{\prime2}}}+f_{1}(w)-d\left\vert w\right\vert
^{-1}w=0, \label{eqe}%
\end{equation}
where $f_{1}$ satisfies (\ref{hfun}); in particular we are interessed by the
case $f_{1}(s)=s.\medskip$

Here  the problem is variational: if $S_{1}(w)$ is
any primitive of $w\mapsto\left\vert w\right\vert ^{b-1}f_{1}(w)$ and
$R(w)=\left\vert w\right\vert ^{b}/b$ if $b\neq0,$ $R(w)=\ln\left\vert
w\right\vert $ if $b=0,$ then (\ref{eqe}) is the Euler equation of the
functional
\[
\mathcal{H(}w,w^{\prime})=\left\vert w\right\vert ^{b-1}\sqrt{w^{2}%
+w^{\prime2}}-S_{1}(w)+dR(w).
\]
Thus the following Painlev\'{e} first integral is constant along the trajectories
\begin{equation}
\CP(w,w')=\frac{\left\vert w\right\vert ^{b+1}}{\sqrt{w^{2}+w^{\prime2}}}-S_{1}%
(w)+dR(w). \label{pai}%
\end{equation}
The system (\ref{FG}) reads as
$$\left\{\BA {l}
w^{\prime}=y\\[2mm]
y^{\prime}=G(w,y)=\myfrac{ bw^{3}%
+(b+1)w\,y^{2}-(f(w)-d|w|^{-1}w)(w^{2}+y^{2})^{3/2}}{w^{2}},\EA\right.
$$
and it is singular on the line $w=0.$ For $w>0$ system (\ref{suw}) reduces  to
\begin{equation}\label{sui}\left\{\BA {l}
w^{\prime}=w\varphi(u)=w\frac{u}{\sqrt{1-u^{2}}}\\[2mm]
u^{\prime}%
=b\sqrt{1-u^{2}}-f_{1}(w)+d.\EA\right.
\end{equation}
In the case $f(w)=w$, the equation satisfied by $u$ is
\begin{equation}
u^{\prime\prime}=(1-b)\frac{u}{\sqrt{1-u^{2}}}u^{\prime}-bu-d\frac{u}%
{\sqrt{1-u^{2}}}. \label{sec}%
\end{equation}

\subsection{Existence of periodic solutions}

From the Painlev\'{e} integral (\ref{pai}), we can describe the solutions, in the phase plane $(w,y).$ Since a complete description is rather long, we
reduce it to the research of periodic solutions.\medskip

\bprop{sig}Let $p=1,$ and consider equation (\ref{eqe}).\smallskip

\noindent (i) If $d\neq0,$ there is no periodic sign changing solution. If $d=0$ there
exists such a solution if and only if $b>-1,$ and then it is unique (up to a translation).\smallskip

\noindent (ii) There exists periodic positive solutions if and only if $b+d>0.$\smallskip

\noindent (iii) Suppose moreover that $f_{1}(w)=w.$ Then the sign changing solution is
given by
\[
w(\tau)=(b+1)\cos(\tau-\tau_{1});
\]
it has period $2\pi$. The orbits $\mathcal{T}_{\left[  (\mu,0)\right]  }$
of the periodic solutions intersect the axis $y=0$ at a first point $(\mu,0)$
such that $\mu<a=b+d,$ and $\mu$ describes $\mu\in\left(  \bar{\mu},a\right)
$ with $\bar{\mu}=0$ if $d\leq0,$ and $\bar{\mu}>0$ if $d>0;$ it is given by
(\ref{mu}),(\ref{H}),(\ref{M}).
\es
\Proof
By symmetry we reduce the study to the case $w\geq 0$ and the painlev\'e integral
(\ref{pai}) takes the form
\begin{equation}
w^{b}\sqrt{1-u^{2}}-S_{1}(w)+dR(w)=C, \label{sit}%
\end{equation}
where we denote 
$$S_{1}(w)=
{\displaystyle\int\nolimits_{0}^{w}}
s^{b-1}f_{1}(s)ds\quad\text{ if  }b>-1,$$
$$S_{1}(w)=
{\displaystyle\int\nolimits_{1}^{w}}
s^{b-1}f_{1}(s)ds+\myfrac{1}{b+1}\quad\text{ if  }b<-1,$$ 
and 
$$S_{1}(w)=
{\displaystyle\int\nolimits_{1}^{w}}
s^{-2}f_{1}(s)ds\quad\text{ if  }b=-1.$$

\noindent {\it Step 1. Periodic sign changing solutions.} The curves in the phase plane $(w,y)$ are given, for $w>0$, by
$$\BA {l}
y^{2}=\left(\myfrac{w^{b+1}}{C-dR(w)+S_{1}(w)}\right)^{2}-w^{2}\\[4mm]
\phantom{y^{2}}
=\myfrac{w^{2}\left(
w^{b}+dR(w)-S_{1}(w)-C\right)  (w^{b}-dR(w)+S_{1}(w)+C)}{\left(
S_{1}(w)+C-dR(w)\right)  ^{2}},%
\EA$$
which defines $\pm y$ in function of $w.$ If there exists a sign changing
periodic solution, the trajectory intersects the axis $w=0$ at some point $(0,\ell)$
with $\ell\geq0,$ thus $y$ needs to $\ell$ as $w$ tends to $0.$ From
(\ref{sit}), it is impossible if $b\leq-1.$ Assume $d\neq0;$ if $-1<b,$ then
near $w=0,$ in any case $y^{2}\leq(b^{2}/d^{2}+1)w^{2},$ thus $\ell=0$ and
$w^{\prime}/w$ is bounded, thus the maximal interval of existence is infinite, 
and we reach a contradiction. If $d=0,$ and $C\neq0,$ then $y^{2}%
=-w^{2}(1+o(1)),$ which is impossible. If $d=C=0,$ then $y^{2}=w^{2}%
(w^{2b}/S_{1}^{2}(w)-1)$ ; observing that the function $w\mapsto\chi
(w)=w^{-b}S_{1}(w)$ is increasing from $0$ to $\infty$, the curve intersects the two
axis at $(0,b+1)$ and ($\chi^{-1}(1),0)$ and this corresponds to a closed orbit.
$\medskip$

\noindent {\it Step 2. Existence of periodic positive solutions.} If we look at the intersection points of any trajectory in the phase
plane with the axis $y=0$, we find that they are given by $H(w)=C,$ where
\[
H(w)=w^{b}+dR(w)-S_{1}(w).
\]
Then $H^{\prime}(w)=w^{b-1}(b+d-f_{1}(w)).$ If $b+d\leq0,$ then $H$ is
decreasing, thus there exist no positive periodic solutions. If $b+d>0,$ the
function $H$ is increasing on $(0,a)$ where it reaches a maximum $M,$ and decreasing on
$(a,\infty)$. The stationary point $(a,0)$ with $a=f_{1}^{-1}(b+d)$
corresponds to $C=M.$ If $b>0,$ then $\lim_{w\rightarrow0}H=0,$ while,  if $b\leq0,$
then $\lim_{w\rightarrow0}H=-\infty.$ Equation $H(w)=C$ has two roots
$0<w_{1}<w_{2}$, if and only if $C\in\left(  \max\{\lim_{w\rightarrow0}%
H,\lim_{w\rightarrow\infty}H\right\} ,M).$ Moreover, if there exists a trajectory
going through $\left(  w_{1},0\right)  $ and $\left(  w_{2},0\right)$ and if one denotes
\[
K(w)=dR(w)-S_{1}(w)=H(w)-w^{b},
\]
one has $K(w)<$ $C$ on $(w_{1},w_{2})$, thus $C>M^{\prime}=\max K=K(f_{1}%
^{-1}(d)).$ Conversely, if
\begin{equation}
 \max\left\{\lim_{w\rightarrow0}H,\lim_{w\rightarrow\infty}H,M^{\prime
}\right\}<C<M,\quad M=H(f^{-1}(b+d)),\quad M^{\prime}=K(f^{-1}(d)),
\label{limc}%
\end{equation}
then there exists a closed orbit going through $\left(  w_{1},0\right)$ and
$\left(  w_{2},0\right)  $.$\medskip$

\noindent {\it Step 3. End of the proof.} The sign changing solution is given by $w^{2}+w^{\prime2}=(b+1)^{2}$ and
its trajectory is a circle with center $0$ and radius $b+1;$ for $w>0,$
$w=(b+1)\sqrt{1-u^{2}}=b\sqrt{1-u^{2}}-u^{\prime},$ thus $u^{\prime}%
=-\sqrt{1-u^{2}},$ and $\theta^{\prime}=-1,$ then $w(\tau)=(b+1)\cos
(\theta-\tau_{1}),$ periodic solution with period $2\pi.$ Now consider the positive periodic
solutions. Here $a=b+d,$ and
\begin{equation}
H(w)=\left\{
\begin{array}
[c]{c}%
(1+d/b)w^{b}-w^{b+1}/(b+1)\text{,\qquad if }b\neq0,-1,\\[2mm]
1+d\ln w-w,\qquad\qquad\qquad\quad\text{ if }b=0,\\[2mm]
(1-d)w^{-1}-\ln w,\qquad\qquad\text{ \qquad if }b=-1,
\end{array}
\right.  \label{H}%
\end{equation}%
\begin{equation}
\left\{
\begin{array}
[c]{c}%
M=a^{b+1}/b(b+1),\quad\quad M^{\prime}=(d^{+})^{b+1}/b(b+1)),\text{\qquad if
}b\neq0,-1,\\[2mm]
M=1+d\ln d-d,\quad\quad M^{\prime}=d\ln d-d,\qquad\quad\qquad\text{ if
}b=0.\\[2mm]
M=-1-ln(d-1),\quad\quad M^{\prime}=-1-\ln d)),\text{\qquad}\qquad\text{if
}b=-1,
\end{array}
\right.  \label{M}%
\end{equation}
If $d\leq0,$ thus $b>0$, then any $C\in\left(  0,M\right)  $ corresponds to a
closed orbit, thus for any $\mu\in\left(  0,a_{1}\right),$ one has a closed
orbit passing through $\left(  \mu,0\right),$ of period still denoted by
$T^{+}(\mu).$ If $d>0,$ in any case, any $C\in\left(  M^{\prime},M\right)  $
corresponds to a closed orbit. If $-1<b<0$, then $H$ is increasing
on $(0,a)$ from $-\infty$ to $M<0,$ and then decreasing on $(a,\infty)$ from
$M$ to $-\infty.$ If $b<-1,$ then $M>M^{\prime}>0.$ If $b<-1,$ then $d>1,$ and
$\lim_{w\rightarrow0}H=-\infty,$ $\lim_{w\rightarrow\infty}H=0,$ $0<M^{\prime
}<M$. If $b=0,$ thus $d>0,$ then $\lim_{w\rightarrow0}H=-\infty,$ then any
$C\in\left(  M-1,M\right)  $ corresponds to a closed orbit. Then $H$ is
increasing on $(0,d)$ from $-\infty$ to $M=1+d\ln d-d\geq0,$
(notice that $M=0\Leftrightarrow d=1$) and then decreasing on $(a_{1},\infty)$
from $M$ to
$-\infty$; let $\bar{\mu}\in\left(  0,b+d\right)  $ be defined by
\begin{equation}
H(\bar{\mu})=M^{\prime}, \label{mu}%
\end{equation}
thus for any $\mu\in\left(  \bar{\mu},a\right),$ one has a closed orbit
passing through $\left(  \mu,0\right),$ with a period still denoted by
$T^{+}(\mu).$ If $-1\leq b<0,$ thus $d>-b>0,$ then $\lim_{w\rightarrow
0}H=-\infty=$ $\lim_{w\rightarrow\infty}H,M<0$, and any $C\in\left(
M^{\prime},M\right)  $ corresponds to a closed orbit (if $b=-1,$ then
$H(w)=(1-d)w^{-1}-\ln w,$ $M=-1-ln(d-1),$ $M^{\prime}=-1-\ln d).$ Returning to
equation (\ref{2dim2}), the conclusion follows with $\bar{\mu}_{q}=\bar{\mu}%
^{q}.$
\qeda

\subsection{Period of the solutions}

Let $p=1,b+d>0.$ Consider the equation (\ref{eqe}). Let $T^{+}(\mu)$ be the
least period of the periodic positive solutions corresponding to the orbit $\mathcal{T}%
_{\left[  (\mu,0)\right]  }$. As in the case $p>1,$ we have a general result:
\begin{equation}
\lim_{\mu\rightarrow a}T^{+}(\mu)=\frac{2\pi}{\sqrt{af_{1}^{\prime}(a)}}.\label{lim-a}
\end{equation}
Next we study the variations of the period in the case of a power $f_{1}(w)=w$.
\bth
{per}Assume $p=1,b+d>0$ and $f_{1}(w)=w.$ Then $\lim_{\mu\rightarrow
a}T^{+}(\mu)=2\pi/\sqrt{b+d}.$ If $d<0,$ then $\lim_{\mu\rightarrow0}T^{+}%
(\mu)=\infty.$ If $d\geq0,$ then $\lim_{\mu\rightarrow\bar{\mu}}T^{+}%
(\mu)=\bar{T}^{+}$ is finite, and given by (\ref{tio}) if $b\notin\{0,-1\},$
by (\ref{bri}) if $b=0,$ and by (\ref{bun}) if $b=-1$. If $d=0,$ then $\bar{T}%
^{+}=\pi(1+1/b).$
\es
\Proof 
{\it Step 1. Assume $b\notin\{0,-1\}$.} From (\ref{sit}), the solutions of
(\ref{eqe}) satisfy
\[
w^{b}\sqrt{1-u^{2}}-\frac{w^{b+1}}{b+1}+\frac{d}{b}w^{b}=C,
\]
thus
\[
u^{\prime}=b\sqrt{1-u^{2}}-w+d=-\sqrt{1-u^{2}}-\frac{d}{b}+\frac{C(1+b)}%
{w^{b}}.%
\]
Eliminating $w$ between the two relations, we find that $Cb(b+1)>0$ and
\[
\left(  d+b\sqrt{1-u^{2}}-u^{\prime}\right)  ^{b/(b+1)}\left(  d+b\sqrt
{1-u^{2}}+bu^{\prime}\right)  ^{1/(b+1)}=(Cb(1+b))^{1/(b+1)}
:=A.\]
When a solution goes through the half-part of its trajectory $\mathcal{T}$ located in $\mathcal{Q}$,
the associated function $u$ increases from $0$ to some $u^{\ast}\in\left(  0,1\right)  $
where the derivative $u^{\prime}$ vanishes and $d+b\sqrt{1-u^{\ast2}}>0;$ next
$d+b\sqrt{1-u^{2}}$ is monotone and positive at $0$ and $u^{\ast},$ thus
$d+b\sqrt{1-u^{2}}>0$ everywhere. And $A=d+b\sqrt{1-u^{\ast2}}=w^{\ast}$
(the value of $w$ when $u=u^{\ast})$. Let
\[
z=\frac{u^{\prime}}{d+b\sqrt{1-u^{2}}}\quad\text{and}\quad G(s)=\left(
1-s\right)  ^{b/(b+1)}\left(  1+bs\right)  ^{1/(b+1)}.
\]
If $b>0,$ then $z\in\left(  -1/b,1\right);$ if $b<-1$ then $z\in\left(
-\infty,1\right);$ if $-1<b<0$ then $z\in\left(  -\infty,1/\left\vert
b\right\vert \right),$ and
\[
G(z)=\frac{A}{d+b\sqrt{1-u^{2}}}.%
\]
Since
\[
G^{\prime}(s)=-bs\left(  1-s\right)  ^{-1/(b+1)}\left(  1+bs\right)
^{-b/(b+1)},\]
and
$$G^{\prime\prime}(s)=-b\left(  1-s\right)  ^{-(b+2)/(b+1)}%
\left(  1+bs\right)  ^{-(2b+1)/(b+1)},
$$
it follows $G(0)=1$, and $0$ is a maximum if $b>0$ and a minimum if $b<0$: if $b>0$,
$G$ increases on $\left(  -1/b,0\right)  $ from $0$ to $1$ and decreases on
$\left(  0,1\right)  $ from 1 to $0;$ if $b<0,$ $G$ decreases on $\left(
-\infty,0\right)  $ from $\infty$ to $1$ and increases on $\left(
0,\min(1,1/\left\vert b\right\vert \right)  )$ from $1$ to $\infty.$ Thus it
has two inverse functions $-L_{1}$ and $L_{2}:$ for $b>0,$ $L_{1}$ maps
$\left(  0,1\right)  $ into $\left(  0,1/b\right)  $ and $L_{2}$ maps $\left(
0,1\right)  $ into $\left(  0,1\right)  ;$ for $b<0,$ $L_{1}$ maps $\left(
1,\infty\right)  $ into $\left(  0,\infty\right)  $ and $L_{2}$ maps $\left(
1,\infty\right)  $ into $\left(  0,\min(1,1/\left\vert b\right\vert \right)
).$ Then%
\begin{equation}
T^{+}=T_{1}^{+}+T_{2}^{+},\qquad T_{i}^{+}=%
{\displaystyle\int\nolimits_{0}^{1}}
\psi_{i,u^{\ast}}(\lambda)d\lambda\label{t12},%
\end{equation}
where%
\[
\psi_{i,u^{\ast}}(\lambda)=\frac{2u^{\ast}}{(d+b\sqrt{1-\lambda^{2}u^{\ast2}%
})L_{i}((d+b\sqrt{1-u^{\ast2}})/(d+b\sqrt{1-\lambda^{2}u^{\ast2}}))}.%
\]

$\bullet$\textbf{ }First suppose $d<0$ (thus $b>0);$ then one looks at the
case where $C\rightarrow0$, thus $\sqrt{1-u^{\ast2}}\rightarrow-d/b,$ thus
$u^{\ast}\rightarrow\bar{u}=\sqrt{1-d^{2}/b^{2}}$ . Near $\bar{u},$
\[
\psi_{i,u^{\ast}}(\lambda)\geq\frac{2u^{\ast}}{b(\sqrt{1-\lambda^{2}u^{\ast2}%
}-\sqrt{1-u^{\ast2}})L_{i}(0)}\geq\frac{-d}{bL_{i}(0)(1-\lambda^{2})},%
\]
therefore $T_{i}^{+}$ tends to $\infty.$

$\bullet$ Suppose $d\geq0,b>0.$ 0ne looks at the case where $C\rightarrow
M^{\prime},$ thus $u^{\ast}\rightarrow1.$ There exists a constant $m>0$ such
that $0\leq1-G(s)=G(0)-G(s)\leq m^{2}s^{2}$ on $\left[  -1/b,1\right]  .$
Indeed $G^{\prime}(0)=0$ and $G^{\prime\prime}$ is bounded on $\left[
-1/2b,1/2\right],$ and on $\left[  -1/b,-1/2b\right]  \cup\left[
1/2,1\right]  $ the quotient $(G(0)-G(s))/s^{2}$ is bounded. Thus
$1/L_{i}(\eta)\leq m/\sqrt{1-\eta}$ on $\left[  0,1\right),$ hence taking
$\eta=(d+b\sqrt{1-u^{\ast2}})/(d+b\sqrt{1-\lambda^{2}u^{\ast2}}),$ and
computing $^{+}$%
\[
1-\eta=\frac{b(1-\lambda^{2})u^{\ast2}}{(d+b\sqrt{1-\lambda^{2}u^{\ast2}%
})\left(  \sqrt{1-u^{\ast2}}+\sqrt{1-\lambda^{2}u^{\ast2}}\right)  },%
\]
one finds $\psi_{i,u^{\ast}}(\lambda)\leq4m/\sqrt{b\left(  1-\lambda
^{2}\right)  }.$ From the Lebesgue theorem, as $u^{\ast}\rightarrow1,$ $T^{+}$
tends to the finite limit
\begin{equation}
\bar{T}^{+}=\bar{T}_{1}^{+}+\bar{T}_{2}^{+},\qquad\bar{T}_{i}^{+}=2%
{\displaystyle\int\nolimits_{0}^{1}}
\frac{d\lambda}{(d+b\sqrt{1-\lambda^{2}})L_{i}(d/(d+b\sqrt{1-\lambda^{2}}))}
\label{tio}%
\end{equation}
in particular if $d=0,$ then $L_{1}(0)=1/b,L_{1}(0)=1,$ thus $\bar{T}%
_{1,1}^{+}=\pi$ and $\bar{T}_{1,2}^{+}=\pi/b.$

$\bullet$ Suppose $b<0,$ thus $d>-b>0.$ Then again $C\rightarrow M^{\prime},$
consequently $u^{\ast}\rightarrow1$. The function
\[
u^{\ast}\rightarrow Q(u^{\ast},\lambda)=\eta=\frac{d+b\sqrt{1-u^{\ast2}}%
}{d+b\sqrt{1-\lambda^{2}u^{\ast2}}}=1-\frac{b(1-\lambda^{2})u^{\ast2}}{\left(
d+b\sqrt{1-\lambda^{2}u^{\ast2}}\right)  (\sqrt{1-u^{\ast2}}+\sqrt
{1-\lambda^{2}u^{\ast2}})}%
\]
is increasing on $\left(  0,1\right)  $ from $1$ to $d/(d+b\sqrt{1-\lambda
^{2}})$ and $d/(d+b\sqrt{1-\lambda^{2}})\leq d/(d+b)=\alpha.$ There exists
$m>0$ such that $0\leq G(s)-1\leq m^{2}s^{2}$ on $\left[  -L_{1}(\alpha
),L_{2}(\alpha)\right],$ thus $1/L_{i}(\eta)\leq m/\sqrt{\eta-1}$ on
$\left(  1/d/(d+b)\right]  .$ Thus as above, $\psi_{i,u^{\ast}}(\lambda
)\leq4m/\sqrt{\left\vert b\right\vert \left(  1-\lambda^{2}\right)  },$ and
$T^{+}$ tends to $\bar{T}^{+}$ defined at (\ref{tio}).

\noindent  {\it Step 2. Assume $b=0$.} There exist periodic solutions for any
$C\in\left(  M-1,M\right).$ The solutions are given by
\[
\sqrt{1-u^{2}}+H(w)=\sqrt{1-u^{2}}+d\ln w-w=C
\]
and $u^{\prime}=-w+d,$ thus $u$ is maximal ($=u^{\ast})$ for $w=d:$ therefore
$\sqrt{1-u^{\ast2}}+H(d)=C,$ then
\[
H(d-u^{\prime})=H(d)+\sqrt{1-u^{\ast2}}-\sqrt{1-u^{2}}%
\]
and $H$ has two inverse functions $H_{i}$ from $\left(  -\infty,H(d)\right)  $
into $(0,d)$ and $(d,\infty),$ thus (\ref{t12}) holds with%
\[
\psi_{i,u^{\ast}}(\lambda)=\frac{2u^{\ast}d\lambda}{(d-H_{i}(H(d)+\sqrt
{1-u^{\ast2}}-\sqrt{1-\lambda^{2}u^{\ast2}})}%
\]
and $\xi=H(d)+\sqrt{1-u^{\ast2}}-\sqrt{1-\lambda^{2}u^{\ast2}}=H(d)-k=H(d+h)$
stays in $\left(  M-1,M\right)  =\left(  H(d)-1,H(d)\right),$ and
$H(d+h)-H(d)\geq-m^{2}h^{2}$ for $H(d+h)\in\left(  M-1,M\right),$ thus
$H(d)-\xi=k\leq m^{2}(d-H_{i}(\xi))^{2},$ thus
\[
\psi_{i,u^{\ast}}(\lambda)\leq\frac{2m}{\sqrt{k}}=\frac{2m\left(
\sqrt{1-u^{\ast2}}+\sqrt{1-\lambda^{2}u^{\ast2}}\right)  }{\sqrt{1-\lambda
^{2}}}\leq\frac{4m}{\sqrt{1-\lambda^{2}}}.%
\]
Therefore, as $u^{\ast}\rightarrow1,$ $T^{+}$ tends to the finite limit
\begin{equation}
\bar{T}^{+}=\bar{T}_{1}^{+}+\bar{T}_{2}^{+},\qquad\bar{T}_{i}^{+}=2%
{\displaystyle\int\nolimits_{0}^{1}}
\frac{d\lambda}{(d-H_{i}(H(d)-\sqrt{1-\lambda^{2}})} \label{bri}.%
\end{equation}

\noindent  {\it Step 3. Assume $b=-1$.} In that case $d>1;$ let $B=-(C+1)\in\left(  \ln(d-1),\ln
d\right)  $ then $B\rightarrow\ln d$ and
\[
u^{\prime}+w=d-\sqrt{1-u^{2}}=(B+1)w-w\ln w=H_{B}(w)
\]
where $H_{B}$ is increasing on $\left(  0,e^{B}\right)  $ from $0$ to $e^{B}$
and decreasing on $\left(  e^{B},\infty\right)  $ from $e^{B}$ to $-\infty;$
it has two inverse functions $L_{B,i}$ from $\left(  -\infty,e^{B}\right)  $
into $(0,e^{B})$ and $(e^{B-1},\infty);$ and $w^{\ast}=d-\sqrt{1-u^{\ast2}%
}=e^{B};$ then (\ref{t12}) holds with
\[
\psi_{i,u^{\ast}}(\lambda)=\frac{2u^{\ast}}{d-\sqrt{1-\lambda^{2}u^{\ast2}%
}-L_{B,i}(d-\sqrt{1-\lambda^{2}u^{\ast2}}}=\frac{2u^{\ast}}{\left\vert
H_{B-1}(L_{B,i}(d-\sqrt{1-\lambda^{2}u^{\ast2}}))\right\vert }.
\]
Because $H_{B-1}(e^{B})=0,$ $H_{B-1}(x)-H_{B-1}(e^{B})=H_{B-1}^{\prime}%
(\xi)(x-e^{B})$ and $x$ ranges onto $\left(
H_{B,1}\left(  d-1\right) ,H_{B,2}\left(  d-1\right)  \right):=\left(  x_{1,B},x_{2,B}\right)$, when
$B\rightarrow\ln d,$ $\left(  x_{1,B},x_{2,B}\right)  \rightarrow\left(
x_{1,\ln d},x_{2,\ln d}\right)$, it follows $\left\vert H_{B-1}^{\prime}%
(\xi)\right\vert \geq1/\mu>0$ independent on $B.$ Moreover $H_{B}%
(x)-H_{B}(e^{B})=(1/2)H_{B}^{\prime\prime}(\xi)(x-e^{B})^{2}=-(1/2\xi
)(x-e^{B})^{2}.$ Thus there exists $m>0$ such that%
\[
H_{B}(x)-H_{B}(e^{B})\leq m^{2}(x-e^{B})^{2}\leq m^{2}\mu^{2}H_{B-1}^{2}(x).
\]
Therefore, near $\ln d,$ taking $x=L_{B,i}(d-\sqrt{1-\lambda^{2}u^{\ast2}}),$ one derive
\[
_{\psi_{i,u^{\ast}}}(\lambda)\leq\frac{2}{m\mu\sqrt{d-\sqrt{1-\lambda
^{2}u^{\ast2}}-e^{B}}}=\frac{2}{m\mu\sqrt{\sqrt{1-u^{\ast2}}-\sqrt
{1-\lambda^{2}u^{\ast2}}}}\leq\frac{4}{m\mu\sqrt{1-\lambda^{2}}}.%
\]
Consequently, as $u^{\ast}\rightarrow1,$ $T_{1,i}^{+}$ tends to the finite limit
\begin{equation}
\bar{T}^{+}=\bar{T}_{1}^{+}+\bar{T}_{2}^{+},\qquad\bar{T}_{1,i}^{+}=2%
{\displaystyle\int\nolimits_{0}^{1}}
\frac{d\lambda}{\left\vert H_{\ln d-1}(L_{\ln d,i}(d-\sqrt{1-\lambda^{2}%
}))\right\vert }. \label{bun}%
\end{equation}
\qeda\\

\noindent\Remark
In the case $d=0,b\neq1,$ notice that $T_{1}^{+}$ and $T_{2}^{+}$ converges to
$\pi/\sqrt{b}$ as $\mu$ tends to $b$ (one can verify it by linearizing the
equation in $u)$ and respectively to $\pi$ and $\pi/b$ as $\mu$ tends to Thus
if those functions are monotonous, they vary in opposite senses and it is not
easy to get the sense of variations of their sum $T_{1}^{+}.$ Moreover in the
phase plane $(w,y_{1}),$ as $\mu$ tends to $0,$ on can observe that the
trajectory tends to a limit curve constituted of a segment $\left[
(0,0),(0,b)\right]  $ and half of the unique closed orbit surrounding $(0,0),$
circle of center $0$ and radius $b+1$, which is covered in a time $\pi
$\bigskip

The case $b=1$ is the most interesting for (\ref{eqe}), since it corresponds
to the initial problem (\ref{2dim2}). In that case we improve the results by
showing the monotonicity of the period function:

\bth{peri} Assume $b=1,$ $d>-1.$ When $d=0$ the period function $T^{+}(\mu)$ is
constant, with value $2\pi,$ thus there exists an infinity of positive solutions
$w$ of (\ref{eqe}), which are all $2\pi$-periodic; they are explicitely given by
\begin{equation}
w=\sqrt{1-K^{\ast2}\sin^{2}\tau}-K^{\ast}\cos\tau,\qquad\tau\in\left[
-\pi,\pi\right] ,\quad K^{\ast}\in\left(  0,1\right)  . \label{giv}%
\end{equation}
When $d\neq0,$ then $T^{+}(\mu)$ is strictly monotone; if $d<0$ it decreases
from $\infty$ to $2\pi/\sqrt{1+d};$ if $d>0$ it increases from
\begin{equation}
\bar{T}^{+}=4%
{\displaystyle\int\nolimits_{0}^{1}}
\frac{du}{\sqrt{(d+\sqrt{1-u^{2}})^{2}-d^{2}}}=4%
{\displaystyle\int\nolimits_{0}^{\pi/2}}
\sqrt{\frac{\cos\theta}{\cos\theta+2d}}d\theta\label{bar}%
\end{equation}
to $2\pi/\sqrt{1+d}.$
\es
\Proof
$\bullet$ If\textbf{ }$d=0$, then $u^{\prime\prime}=-u,$ from (\ref{sec}), and
$u=\sin\theta\in\left[  0,1\right),$ thus the positive solutions $w$ are
given in $\mathcal{Q}$ by
\[
u=K^{\ast}\sin\tau,\qquad K^{\ast}\in\left[  0,1\right) ,\qquad\tau\in\left[
0,\pi\right] ,
\]
and the period $T^{+}$ is constant, equal to $2\pi$. We obtain an infinity of
positive solutions $w$, given explicitely by
\[
w=\sqrt{1-u^{2}}-u^{\prime}=\sqrt{1-K^{\ast2}\sin^{2}\tau}-K^{\ast}\cos
\tau,K^{\ast}\in\left(  0,1\right)
\]
which intersect the axe $y=0$ at points $w_{i}=(1\mp K^{\ast}).$

$\bullet$ In the general case $d>-1,$ we find
\[
\left(  d+\sqrt{1-u^{2}}-u^{\prime}\right)  \left(  \sqrt{1-u^{2}}+u^{\prime
}+d\right)  =A^{2}%
\]
that means $G$ is symmetric: $G(s)=\sqrt{1-s^{2}},$ thus
\[
u^{\prime2}=(d+\sqrt{1-u^{2}})^{2}-A^{2}%
\]
$\sqrt{1-u^{\ast2}}=A-d=\sqrt{2C}-d;$ thus here $T_{1}^{+}=T_{2}^{+},$ and
\[
T^{+}=4%
{\displaystyle\int\nolimits_{0}^{1}}
\frac{d\lambda}{\sqrt{\Psi(u^{\ast2},\lambda)}},\quad\text{where}\quad
\Psi(s,\lambda)=\frac{(d+\sqrt{1-\lambda^{2}s})^{2}-(d+\sqrt{1-s})^{2}}{s}.
\]
We show that the period function is strictly monotone with respect to
$u^{\ast}$. Because %
\[
s^{2}\frac{\partial\Psi}{\partial s}(s,\lambda)=d(d+1)\left(  1/\sqrt
{1-s}-1/\sqrt{1-\lambda^{2}s}\right)  >0,
\]
we see that $T^{+}$ is increasing if $d<0$ and decreasing if $d>0$ (and we find again
that it is constant if $d=0)$. Also $\mu$ can be expressed explicitely in
terms of $u^{\ast}$ by
\[
\mu=d+1-\sqrt{(d+1)^{2}-(d+\sqrt{1-u^{\ast2}})}\,.%
\]
Therefore $\mu$ is decreasing with respect to $u^{\ast},$ hence $T^{+}$ is
decreasing with respect to $\mu$ if $d<0$ and increasing if $d>0.$
\qeda
\subsection{Returning to the initial problem}

{\it Proof of Theorem 3.} Here $\alpha_{q}=\beta_{q}=1/q,$ the
substitution (\ref{cvq}) takes the form $\omega\left(  \sigma\right)
=\left\vert w(\sigma)\right\vert ^{1/q-1}w(\sigma),$ and thus $b=1,$ and $d=c$
from (\ref{b1d1}). Then the existence of sign changing solutions of (\ref{2dim2})
is given by \rprop{sig}. The constant solutions exist whenever
$c+1>0.$ Next we look for positive solutions of smallest period $2\pi/k$ applying \rth{per} and \rth{peri}. If $c<0$ the period function $T^{+}$
decreases from $\infty$ to $2\pi/\sqrt{1+c}>2\pi,$ thus there exist no
solution. If $c>0,$ $T^{+}$ increases from $\bar{T}^{+}$ given by (\ref{bar})
to $2\pi/\sqrt{1+c},$ thus it takes once any intermediate value, which gives
one solution (up to a translation) for any $k\in\left(  k_{1},k_{2}\right)  $.
If $c=0,$ the solutions $\omega_{K}$ are given explicitely by (\ref{giv}), and
$\omega_{0}^{+}$ is obtained from $\omega_{0};$ this means that system (\ref{FG})
does not satisfy the uniqueness property at $(0,0).$
\qeda

\end{document}